\documentclass[11pt, reqno]{amsart}

\usepackage{amsmath,amssymb,amscd,amsthm,amsxtra}
\usepackage[dvips]{graphics,epsfig}

\allowdisplaybreaks[2]

\newtheorem{theorem}{Theorem} [section]

\newtheorem{lemma}[theorem]{Lemma}
\newtheorem{proposition}[theorem]{Proposition}
\newtheorem{remark}[theorem]{Remark} 

\DeclareMathOperator*{\intt}{\int}

\DeclareMathOperator{\MAX}{MAX}

\newcommand{\noi}{\noindent}
\newcommand{\Z}{\mathbb{Z}}
\newcommand{\R}{\mathbb{R}}

\newcommand{\T}{\mathbb{T}}

\newcommand{\al}{\alpha}
\newcommand{\dl}{\delta}
\newcommand{\Dl}{\Delta}
\newcommand{\eps}{\varepsilon}

\newcommand{\g}{\gamma}
\newcommand{\G}{\Gamma}
\newcommand{\ld}{\lambda}

\newcommand{\s}{\sigma}
\newcommand{\ft}{\widehat}
\newcommand{\wt}{\widetilde}
\newcommand{\cj}{\overline}
\newcommand{\dx}{\partial_x}
\newcommand{\dt}{\partial_t}

\newcommand{\jb}[1]
{\langle #1 \rangle}

\numberwithin{equation}{section}
\numberwithin{theorem}{section}

\begin{document}

\title
[nonlinear smoothing under randomization]
{\bf Remarks on nonlinear smoothing under randomization for the periodic KdV
and the cubic Szeg\"o equation}

\author{Tadahiro Oh}

\address{Tadahiro Oh\\
%
%
Department of Mathematics\\
Princeton University\\
Fine Hall, Washington Rd\\
Princeton, NJ 08544-1000, USA}

\email{hirooh@math.princeton.edu}

\thanks{
T.O. would like to thank D\'epartement de Math\'ematiques d'Orsay at  Universit\'e Paris-Sud,
where a part of this manuscript was prepared.}



\begin{abstract}
We consider Cauchy problems of 
some dispersive PDEs with random initial data.
In particular, 
we construct local-in-time solutions
to the mean-zero periodic KdV
almost surely for the initial data in the support of the mean-zero Gaussian
measures on $H^s(\T)$, $s > s_0$
where $s_0 = -\frac{11}{6}+\frac{\sqrt{61}}{6} \thickapprox -0.5316 <-\frac{1}{2}$,
by exhibiting nonlinear smoothing under randomization
on the second iteration of the integration formulation.
We also show that there is no nonlinear smoothing
for the dispersionless cubic Szeg\"o equation
under randomization of initial data.
\end{abstract}

\maketitle

\tableofcontents

\section{Introduction}

\subsection{Nonlinear smoothing under randomization on initial data
for nonlinear Schr\"odinger equations} \label{SUBSEC:1.1}
In studying  invariance of the Gibbs measure for 
the defocusing cubic  nonlinear Schr\"odinger equation (NLS) on $\T^2$, 
Bourgain \cite{BO7} considered the following Cauchy problem
for the 2-$d$ Wick ordered cubic NLS:{\footnote{
In \cite{BO7}, Wick ordering is a renormalization needed for constructing
the Gibbs measure on $\T^2$.
The equation \eqref{NLS0} appears
as an equivalent formulation of the NLS 
obtained from the Wick ordered Hamiltonian.}
\begin{equation} \label{NLS0}
\begin{cases}
i u_t - \Delta u \pm (u |u|^2 -2 u \int \hspace{-9.8pt} - |u|^2 dx)= 0 \\
u|_{t= 0} = u_0,
\end{cases}
\quad x \in \T^2 = \R^2/(2\pi\Z)^2
\end{equation}

\noi
with {\it random} initial data $u_0$ of the form 
\begin{equation} \label{IV0}
u_0(x) = u_0^\omega(x) = \sum_{n \in \mathbb{Z}^2} \frac{g_n(\omega)}{\sqrt{1+ |n|^2}} e^{in\cdot x},
\end{equation}

\noi
where $\{ g_n\}_{n \in \mathbb{Z}^2} $ 
is a family of independent standard complex-valued Gaussian random variables
on a probability space $(\Omega, \mathcal{F}, P)$.
We can regard
$u_0$ in \eqref{IV0} as a typical element
in the support of the Gaussian part 
 of the Gibbs measure on $\T^2$:\footnote{We introduced
$-\frac{1}{2}\int |u|^2 dx$ 
to avoid the zero frequency issue. This modification is done implicitly in \cite{BO7}.}
\begin{equation} \label{Gaussian0}
d \rho = Z^{-1} \exp \Big(-\frac{1}{2}\int |u|^2 dx -\frac{1}{2} \int |\nabla u|^2 dx\Big)
\prod_{x \in \mathbb{T}^2} d u(x).
\end{equation}

\noi
Note that \eqref{Gaussian0} is basically the Wiener measure on $\T^2$.
Hence, 
we see that $u_0$ of the form \eqref{IV0}
belongs  almost surely (a.s.) to $H^{s}(\mathbb{T}^2) \setminus L^2(\mathbb{T}^2)$
for any $s<0$.\footnote{This regularity 
can be easily determined  
from the basic theory of Gaussian measures on Hilbert and Banach spaces.
See Kuo \cite{KUO} and Zhidkov \cite{Z}.}

In \cite{BO1}, Bourgain 
introduced a new weighted space-time Sobolev space
$X^{s, b}(\mathbb{T}^d \times \mathbb{R})$,
whose norm is given by 
\[\| u\|_{X^{s, b}(\mathbb{T}^d \times \mathbb{R})} 
= \|\jb{n}^s \jb{\tau - |n|^2}^b \ft{u}(n, \tau)\|_{l^2_nL^2_\tau  (\mathbb{Z}^d \times\mathbb{R})}, \]

\noi
where $\jb{ \: \cdot \:} = 1 + |  \cdot  | $, and proved that \eqref{NLS0} is locally well-posed
in $H^s(\T^2)$ for $ s> 0$.
This result barely misses the scaling-critical space $L^2(\T^2)$
for \eqref{NLS0} due to a slight loss of derivative
in the periodic $L^4$-Strichartz estimate on $\T^2$.

There are two main components in establishing  invariance of Gibbs measures
 for Hamiltonian PDEs:
  \begin{itemize}
  \item[(a)] construction of the Gibbs measure,
 \item[(b)] (at least local-in-time) well-posedness on the support of the Gibbs measure.
 \end{itemize}
 
\noi 
 In \cite{BO7}, Bourgain constructed the Gibbs measure
after applying Wick ordering to the nonlinear part of the Hamiltonian
(in the defocusing case.)
As for (b), 
one needs to construct 
a continuous flow of \eqref{NLS0} on the support of the Gibbs measure, 
i.e.~outside $L^2(\T^2)$. 
This is exactly the main difficulty in \cite{BO7},
due to the absence of deterministic well-posedness even in $L^2(\T^2)$. 
In order to resolve this issue, 
Bourgain considered the Cauchy problem
\eqref{NLS0} with random initial data \eqref{IV0}
and successfully constructed local-in-time solutions almost surely in $\omega$
by exhibiting {\it nonlinear smoothing under randomization of initial data}. 
Then, by invariance of the finite dimensional  Gibbs measures
with approximation argument,
such local-in-time solutions were then extended globally in time.

We briefly describe Bourgain's idea in the following.
First, write \eqref{NLS0} in the Duhamel formulation:
\begin{equation} \label{NLS3}
u(t) = \G u(t) := S(t) u_0 \pm i \int_0^t S(t - t')  \mathcal{N}(u) (t') d t',
\end{equation}

\noi
where
$S(t) = e^{-i t \Delta  }$, $u_0$ is as in \eqref{IV0},  and $\mathcal{N}(u)  = u |u|^2 - 2u \int \hspace{-9.8pt} - |u|^2$.
Note that the linear part $S(t) u_0$ has the same regularity as $u_0$ for each fixed $t \in \mathbb{R}$,
i.e.
$S(t) u^\omega_0 \notin L^2 (\mathbb{T}^2)$ a.s.
However, by a combination of deterministic PDE theory 
and probabilistic techniques, 
Bourgain showed that 
$\int_0^t S(t - t')  \mathcal{N}(u) (t') d t'$ lies almost surely
in a smoother space $H^s (\mathbb{T})$ for some small $s  > 0$
(namely,  nonlinear smoothing under randomization.)
Indeed, he showed that 
for each small $T> 0$, 
there exists a set $\Omega_T$ with complemental measure $< e^{-\frac{1}{T^\dl}}$
such that  
for $\omega \in \Omega_T$,
$\G$ defined in \eqref{NLS3} is a contraction
on a ball around the linear solution,
i.e. on $S(t) u_0^\omega + B$,
where $B$ denotes the ball of radius 1 in 
the usual Bourgain space $X^{s, \frac{1}{2}+, T}$
for some small $s  > 0$,
where $X^{s, \frac{1}{2}+, T}$ is a local-in-time version of 
the $X^{s, b}$ space on $[-T, T]$.

Recently, several results appeared in this direction,
exhibiting  nonlinear smoothing under randomization of initial data.
See Burq-Tzvetkov \cite{BT2} for the cubic nonlinear wave equation on a three-dimensional compact
Riemannian manifold and Thomann \cite{LT} for NLS with a confining potential on $\R^d$.
(In \cite{LT}, there is a statement for NLS without a potential
but the result is stated in terms of the Sobolev space corresponding to the Laplace operator 
with a confining potential.)
Colliander-Oh \cite{CO1} considered
the 1-$d$ Wick ordered cubic NLS \eqref{NLS0} on $\T$
(both defocusing and focusing)
with random initial data $u_0$ of the form:
\begin{equation} \label{IV}
u_0(x) = u_0^\omega(x) = \sum_{n \in \mathbb{Z}} \frac{g_n(\omega)}{\sqrt{1+ |n|^{2\al}}} e^{inx}.
\end{equation}

\noi
Note that $u_0$ in \eqref{IV} is a.s.~in $H^{\al-\frac{1}{2}-}(\T):=\bigcap_{s < \al - \frac{1}{2}} H^s(\T) 
\setminus H^{\al - \frac{1}{2} }(\T)$.
In particular, $u_0$ is almost surely in the negative Sobolev spaces for $\al \leq \frac{1}{2}$.
Also, recall that $u_0$ in \eqref{IV} represents a typical element
in the support of the Gaussian measure $\rho_\alpha$ on the distributions on $\T$:
\begin{equation} \label{Gaussian1}
d \rho_\al = Z^{-1} \exp \Big(-\frac{1}{2}\int |u|^2 dx -\frac{1}{2} \int |D^{\al} u|^2 dx\Big)
\prod_{x \in \mathbb{T}} d u(x),
\end{equation}

\noi
where $D = \sqrt{-\dx^2}$.
In \cite{CO1}, it is shown that \eqref{NLS0}
is locally well-posed almost surely in $H^{\al-\frac{1}{2}-}(\T)$ for each $\al > \frac{1}{6}$,
i.e.~in $H^s(\T)$ for each $s > -\frac{1}{3}$.
Moreover,  we constructed 
almost sure global-in-time solutions \eqref{NLS0}
for each $\al > \frac{5}{12}$, i.e.~in $H^s(\T)$ for each $s > -\frac{1}{12}$.

The local-in-time argument in \cite{CO1}
closely follows that of Bourgain \cite{BO7}.
The main goal is to show that 
for each small $T>0$, there exists a set $\Omega_T$ 
with $P(\Omega_T^c)<e^{-\frac{1}{T^\dl}}$
such that
for $\omega \in \Omega_T$,
$\G$ in \eqref{NLS3} is a contraction
on $S(t) u_0^\omega + B$, 
where $B$ denotes the ball of radius 1 in 
$X^{s, \frac{1}{2}+, T}$
for some $s  \geq 0$.
By the standard linear estimate on the Duhamel term in \eqref{NLS3},  
it suffices to prove
\begin{equation} \label{trilinear0}
\| \mathcal{N}(u)\|_{X^{s, -\frac{1}{2}+, T}} \lesssim T^\theta, 
\end{equation}

\noi
for $u \in S(t) u_0^\omega + B$
with some $\theta > 0$,
where 
$S(t) = e^{-i t \dx^2  }$, $u_0^\omega$ is as in \eqref{IV},
and 
$\mathcal{N}(u) = \mathcal{N}(u, u, u) = u |u|^2 - 2u \int \hspace{-9.8pt} - |u|^2$.
Since $u \in S(t) u_0^\omega + B$,
we can write $u$ as 
$u = S(t) u_0^\omega + v$ for some $v$ with $ \|v \|_{X^{s, \frac{1}{2}+, T}} \leq 1$.
Hence, it suffices to show  \eqref{trilinear0}
for $\mathcal{N} (u_1, u_2, u_3)$, 
assuming that  $u_j$ is either of the type
\begin{itemize}
\item[(I)] linear part: {\it random, less regular}
\[u_j (x, t) =  S(t) u_0^\omega
=  \sum_{n \in \mathbb{Z} } \frac{ g_n(\omega)}{\sqrt{1+|n|^{2\al}}} e^{i(nx + n^2t)}, \ \ \text{or}
\]
\item[(II)] nonlinear part: {\it deterministic, smoother}
\[u_j = v_j \text{ with } \|v_j \|_{X^{s, \frac{1}{2}+, T}} \leq 1.\]
\end{itemize}

\noi
Then, \eqref{trilinear0}
was established by a combination of deterministic
multilinear analysis and probabilistic argument;
in particular, the hypercontractivity of the Ornstein-Uhlenbeck semigroup 
related to products of Gaussian random variables
played an important role.

\bigskip

In the following subsections, we discuss our main results.
The first result is on KdV equation
with a derivative nonlinearity.
We show that 
a simple application of the ideas in probabilistic Cauchy theory
\cite{BO7, CO1} with a fixed point argument fails (Subsection \ref{SUBSEC:4.1}.)
Nonetheless, we apply the {\it second iteration argument}\footnote{The second iteration argument
is a name of the method, where we (partially) iterate the Duhamel formulation 
as in \eqref{KDVduhamel3} and \eqref{KDVduhamel4}. This is not to be confused with the {\it second iterate} as in \eqref{Sz2} and \eqref{KDV2}.}
in the probabilistic setting to improve known deterministic well-posedness results
(without using complete integrability), i.e. $s = -\frac{1}{2}$.
See Theorem \ref{THM:LWP}.
The second result is on the {\it dispersionless} cubic Szeg\"o equation (See \eqref{Szego} below.)
Here, we show that 
even with randomization on initial data, one can not improve the deterministic well-posedness result
(Proposition \ref{PROP:Szego}.)
This result in particular indicates that 
it is important to have both randomization 
and dispersion together to yield an improvement over deterministic results.

\subsection{Korteweg-de Vries equation} 
\label{SUBSEC:1.2}
In this part, 
we  consider the periodic Korteweg-de Vries (KdV) equation:
\begin{equation} \label{KDV}
\begin{cases}
u_t + u_{xxx} +  u u_x  = 0 \\ 
u \big|_{t = 0} = u_0,
\end{cases}
\quad x \in \T = \R/2\pi \Z
\end{equation}

\noi
with random initial data $u_0$ of the form
\begin{equation} \label{IVKDV}
u_0(x) = u_0^\omega(x) = \sum_{n\ne 0} \frac{g_n(\omega)}{|n|^\al} e^{inx}
\in H_0^{\al - \frac{1}{2}-}(\T) \ \ \text{a.s.},
\end{equation}

\noi
where $\{g_n\}_{n \in\mathbb{N}}$ is a family of independent
standard complex-valued Gaussian random variables
with $g_{-n} = \cj{g_{n}}$ for $n \in \mathbb{N}$
and $H^s_0(\T)$ denotes a subspace of the Sobolev space $H^s(\T)$
consisting of real-valued mean-zero elements.
Note that  $u_0$ in \eqref{IVKDV} represents a typical element
in the support of the Gaussian measure on the real-valued mean-zero distributions on $\T$:
\begin{equation} \label{Gaussian2}
d \rho_{0, \al} = Z^{-1} \exp \Big( -\frac{1}{2} \int (D^{\al} u)^2 dx\Big)
\prod_{x \in \mathbb{T}} d u(x), \quad u, \text{ mean 0}.
\end{equation}

\noi
Our main goal is to construct local-in-time solutions
in $C([-T, T];H^{\al-\frac{1}{2}-}_0(\T))$
for each $\al \in (\al_0, 0]$ with some $\al_0 < 0$
by exhibiting nonlinear smoothing under randomization on initial data.

\medskip

First, let us briefly review recent well-posedness results on KdV \eqref{KDV}.
By introducing the $X^{s, b}$ space adapted to (the linear part of) KdV 
with the norm 
\begin{equation} \label{Xsb}
\| u \|_{X^{s, b}(\mathbb{T} \times \mathbb{R})} = \| \jb{n}^s \jb{\tau - n^3}^b 
\ft{u}(n, \tau) \|_{l^2_n L^2_\tau(\mathbb{Z} \times \R)},
\end{equation}

\noi
Bourgain \cite{BO11} proved local well-posedness of \eqref{KDV} in $L^2(\mathbb{T})$
via the fixed point argument,
immediately yielding  global well-posedness in $L^2(\mathbb{T})$
thanks to the conservation of the $L^2$-norm.
Kenig-Ponce-Vega \cite{KPV4} (also see \cite{CKSTT4}) improved Bourgain's result 
and proved local well-posedness in $H^{-\frac{1}{2}}(\T)$
by establishing the bilinear estimate 
\begin{equation} \label{KPVbilinear}
\| \dx(uv) \|_{X^{s, -\frac{1}{2}}} \lesssim \| u \|_{X^{s, \frac{1}{2}}} \| v \|_{X^{s, \frac{1}{2}}}, 
\end{equation}

\noi
for $s \geq -\frac{1}{2}$ under the mean zero assumption on $u$ and $v$.
Colliander-Keel-Staffilani-Takaoka-Tao \cite{CKSTT4} proved 
the corresponding global well-posedness result via the $I$-method. 

There are also results on \eqref{KDV} which exploit its complete integrability.
In \cite{BO3}, Bourgain proved global well-posedness of \eqref{KDV}
in the class $M(\T)$ of measures $\mu$, assuming that 
its total variation $\|\mu\|$ is sufficiently small.
His proof is based on the trilinear estimate for the {\it second iteration} of 
the integral formulation of \eqref{KDV}, 
assuming an a priori uniform bound on the Fourier coefficients of the solution $u$ of the form 
\begin{equation} \label{BOO}
\sup_{n\in \mathbb{Z}} |\ft{u}(n, t)| < C
\end{equation}

\noi
for all $t\in \R$. 
Then, he established \eqref{BOO} using the complete integrability.
More recently, Kappeler-Topalov \cite{KT} proved global well-posedness of the KdV in $H^{-1}(\T)$
via the inverse spectral method.

Next, we state results regarding the necessary conditions\footnote{These results
are often referred to as ill-posedness results.} on the  regularity 
with respect to smoothness or uniform continuity 
of the solution map $: u_0 \in H^s (\mathbb{T}) \to u(t) \in H^s(\mathbb{T})$.
Bourgain \cite{BO3} showed that if the solution map is $C^3$, 
then $s \geq -\frac{1}{2}$.  
Christ-Colliander-Tao \cite{CCT}
proved that if the solution map is uniformly continuous, 
then $s \geq -\frac{1}{2}$.
(Also, see Kenig-Ponce-Vega \cite{KPV5}.) 
These results, in particular, state that
it is non-trivial to construct solutions in 
$ C([-T, T];H^{s})$ for  $s  <-\frac{1}{2}$.


\medskip

In the following, 
we combine deterministic and probabilistic arguments
to
construct local-in-time solutions 
$u \in C([-T, T];H^{s})$, $s = \al - \frac{1}{2}- <-\frac{1}{2}$,
with initial data \eqref{IVKDV}, where $T = T(\omega)$.
The main difficulty is to lower the differentiability
so that 
$\al \leq 0$, i.e. $s < -\frac{1}{2}$. 
Our focus is to develop a method
{\it without} complete integrability of \eqref{KDV},
exploiting nonlinear smoothing under randomization.
It turns out that when $\al \leq 0$, 
Bourgain's idea in \cite{BO7} - 
a fixed point argument around the linear solution
discussed in Subsection \ref{SUBSEC:1.1}
(as in \cite{BO7, BT2, CO1, LT}) -  
does {\it not} work for \eqref{KDV}.
See Subsection \ref{SUBSEC:4.1}.
Instead, we adapt the {\it second iteration} argument
in the probabilistic setting.

\begin{theorem} \label{THM:LWP}
Let $a_0$ be the positive root of $a^2 + \frac{5}{3}a - 1 = 0$
given by
\begin{equation} \label{a0}
a_0 = -\tfrac{5}{6} +\tfrac{\sqrt{61}}{6} \thickapprox 0.4684,
\end{equation}

\noi
and let $\al_0 := a_0 - \frac{1}{2} \thickapprox -0.0316$.
Then, 
for each $\al \in (\al_0, 0]
\thickapprox (-0.0316, 0]$,
KdV \eqref{KDV} is locally well-posed
almost surely in $H_0^{\al-\frac{1}{2}-}(\T)$.
More precisely, there exist $c, \beta > 0$ such that 
for each $T \ll 1$, there exists a set $\Omega_T \in \mathcal{F}$ with the following properties:

\begin{enumerate}
\item[(i)] $P(\Omega_T^c) = \rho_{0,\al} \circ u_0(\Omega_T^c) < e^{-\frac{c}{T^\beta}}$,
where $u_0:\Omega \to H_0^{\al-\frac{1}{2}-}(\mathbb{T})$.

\item[(ii)] For each $\omega \in \Omega_T$
there exists a (unique) solution $u$ of \eqref{KDV} in 
\[e^{- t \dx^3 }u_0 + C([-T, T];H_0^{-\frac{1}{2}+}(\mathbb{T})) \subset C([-T, T];H_0^{\al - \frac{1}{2}-}(\mathbb{T}))\]
with the initial condition $u_0^\omega$
given by \eqref{IVKDV}.
Here, the uniqueness holds only in a mild sense.
See Subsection \ref{SUBSEC:4.3}.
\end{enumerate}

\noi
In particular, we have almost sure local well-posedness of KdV with respect to the Gaussian measure
$\rho_{0, \al}$ in \eqref{Gaussian2} supported on $H^s(\T)$ for each $s >s_0$
with 
\begin{equation} \label{s0}
s_0 = -\tfrac{11}{6}+\tfrac{\sqrt{61}}{6}\thickapprox -0.5316 .
\end{equation}

\end{theorem}

%
%

\noi
The novelty of Theorem \ref{THM:LWP} 
is the construction of local-in-time solutions  in 
$C([-T, T];H_0^{s}(\mathbb{T}))$ with $s < -\frac{1}{2}$
 without using complete integrability.
This is non-trivial in view of ill-posedness results in $H^s(\T)$, $s <-\frac{1}{2}$,
described above.
The argument combines the second iteration 
and the probabilistic argument, 
and does not rely on complete integrability
of the equation. Therefore, it can be apply to other non-integrable KdV variants.
Lastly, the basic argument in the previous work such as \cite{BO7, BT2, CO1, LT} in random Cauchy theory
is based on 
a fixed point argument around the (random) linear solution
discussed in Subsection \ref{SUBSEC:1.1}.
(This can be also viewed as a fixed point argument for the nonlinear part,
regarding the linear part as a random ``forcing term''.)
As shown in Subsection \ref{SUBSEC:4.1}, 
an attempt to follow this argument fails for KdV.
The proof of Theorem \ref{THM:LWP} presents 
a new way (via non-``fixed point'' argument - namely the second iteration argument
in the probabilistic setting) to construct solution in random Cauchy theory.
See Richards \cite{R} 
for a successful application of this idea to the random Cauchy theory of the quartic KdV
(with nonlinearity $u^3 u_x$.)

%
%
%


\begin{remark} \rm
The regularity $s_0 \thickapprox -0.5316$ in Theorem \ref{THM:LWP} is by no means sharp.
One of the main purposes of this result is to show what one can do when Bourgain's idea 
on random Cauchy theory in \cite{BO7}
 - a fixed point argument around the linear solution - fails.
In Section \ref{SEC:4}, we adapt Bourgain's deterministic second iteration argument
\cite{BO3} in the probabilistic setting.
In particular, 
we use the (deterministic) estimates in \cite{BO3}, that are valid for $s<-\frac{1}{2}$,
without modification.
However,  the values of $s$ and $b$ of the $X^{s, b}$-norm in these estimates
from \cite{BO3}
are strongly interrelated,
giving a restriction the regularity.
See Section \ref{SEC:4}.
We believe that one can lower the regularity $s_0$ in Theorem \ref{THM:LWP}
by improving the deterministic estimates in \cite{BO3}
with the values of $s$ and $b$ ``independent''.
(Namely, keep $b = \frac{1}{2}-$ or $\frac{1}{2}+$
and determine the value of $s$ for which the estimates hold.)
We do not pursue this direction. 
Lastly, note that Theorem \ref{THM:LWP}
covers the critical value $\al = 0$.  See Remark \ref{REM:al}.
\end{remark}

\begin{remark} \label{REM:al}\rm
The value $\al = 0$, corresponding to the regularity $s = -\frac{1}{2}-$, 
is the critical value
in terms of the know well-posedness/ill-posedness results discussed above.
Moreover, 
when $\al = 0$, $u_0$ in \eqref{IVKDV} corresponds to the (mean-zero) Gaussian white noise on $\T$,
which is formally invariant in view of the $L^2$-conservation.
Then, one can follow Bourgain's argument in \cite{BO7}
(i) to extend the solutions in Theorem \ref{THM:LWP}
globally in time, i.e.~in $C(\R_t; H^{-\frac{1}{2}-}(\T))$, and 
(ii) to establish invariance of the white noise under the flow of KdV.
See Quastel-Valk\'o \cite{QV}, 
Oh \cite{OH4, OH7}, and 
Oh-Quastel-Valk\'o \cite{OQV} for other proofs of invariance of white noise for KdV.
Recently, Richards \cite{R} 
applied the second iteration argument to the quartic KdV
 in the probabilistic setting
and constructed solutions below the deterministic threshold $s = \frac{1}{2}$
(and established invariance of the Gibbs measure for the quartic KdV.)
\end{remark}

\begin{remark} \rm

In \cite{OH4, OH6, OH7}, we proved local well-posedness of \eqref{KDV}
in $\mathcal{F} L^{s, p}(\T)$ with $s > -\frac{1}{2}, p = 2+$, $sp < -1$,
where the Fourier-Lebesgue space $\mathcal{F} L^{s, p}(\T)$
is defined by the norm
\begin{equation} \label{FLP}
 \|f\|_{\mathcal{F} L^{s, p}(\T)} = \| \jb{n}^s \ft{f}(n)\|_{l^p_n(\Z)}.
 \end{equation}

\noi
The proof is based on the second iteration argument introduced in \cite{BO3}.
It is known \cite{Benyi, OH4} that $u_0$ in \eqref{IVKDV} is almost surely in $\mathcal{F} L^{s, p}$
for $(s-\al) p < -1$.
Hence, this result shows that there exists $\al_0 < 0$
such that for each $\al \in [\al_0, 0]$, 
KdV \eqref{KDV} is a.s.~locally well-posed
with initial data \eqref{IVKDV}
such that $u \in C([-T, T]; \mathcal{F} L^{s, p})$ with $s = \al - \frac{1}{p}- >-\frac{1}{2}$,
where $T = T(\|u_0(\omega)\|_{\mathcal{F} L^{s, p}})$.
We point out that this deterministic well-posedness argument in $\mathcal{F}L^{s, p}$ completely fails when $p = 2$
and $s <-\frac{1}{2}$.
In particular, Theorem \ref{THM:LWP}
does not follow from the result in \cite{OH4, OH6, OH7}.

\end{remark}

\begin{remark}\rm
A linear part of a solution constructed in Theorem \ref{THM:LWP}
indeed lies in $C([-T, T];B(\mathbb{T}))$
for any Banach space $B (\T) \supset H_0^\al(\T)$
such that $(H_0^\al, B, \rho_{0, \al})$ is an abstract Wiener space
(roughly speaking, any Banach space $B$ containing $H_0^\al$ where the Gaussian measure $\rho_{0, \al}$
makes sense as a countable additive probability measure.)\footnote{Strictly speaking, we need
that $e^{- t \dx^3 }$ acts on $B (\T)$ continuously.
e.g. we need $p <\infty$ for $\mathcal{F}L^{s, p}$.}
In this case, a solution $u$ to \eqref{KDV}
lies in 
\[ u = e^{- t \dx^3 }u_0 + (-\dt - \dx^3)^{-1} u \in
C([-T, T];B(\mathbb{T})) + C([-T, T];H_0^{-\frac{1}{2}+}(\mathbb{T})).\]

As examples of $B$, 
we can take the Sobolev spaces $W^{s, p}$ with $ s < \al - \frac{1}{2}$, and 
the Fourier-Lebesgue spaces $\mathcal{F}L^{s, p}$ with $s < \al - \frac{1}{p}$,
where $\mathcal{F}L^{s, p}$ is defined via the norm \eqref{FLP}.
See B\'enyi-Oh \cite{Benyi}
for regularity of $\rho_{0, \al}$ in different function spaces.
We can also take the Besov spaces $B^{\al - \frac{1}{2}}_{p, \infty}$ with $p < \infty$.
In \cite{Benyi}, we study the regularity of $\rho_\al$ in \eqref{Gaussian1} with $\al = 1$
but it can be easily adjusted for $\rho_{0, \al}$ in \eqref{Gaussian2} with any $\al$.

\end{remark}

\subsection{Cubic Szeg\"o equation} \label{SUBSEC:1.3}

In studying the cubic NLS: $i u_t - \Dl u = |u|^2 u$ on a manifold $M$,
Burq-G\'erard-Tzvetkov \cite{BGT} observed that dispersive properties are strongly influenced
by the geometry of the underlying manifold $M$. 
G\'erard-Grellier \cite{GG} further pointed out 
that dispersion disappears completely when $M$ is a sub-Riemannian manifold.
As a toy model to study {\it non-dispersive} Hamiltonian equations,
G\'erard-Grellier \cite{GG} introduced the {\it cubic Szeg\"o equation}:
\begin{equation} \label{Szego}
\begin{cases}
i u_t = \Pi (|u|^2u)\\
u|_{t= 0} = u_0
\end{cases}
\quad x \in \T,
\end{equation}

\noi
where $\Pi$ is the Szeg\"o projector
onto the non-negative frequencies,
i.e. \[\Pi (f) := \sum_{n\geq0} \ft{f}(n) e^{inx}.\]

\noi
It turned out that \eqref{Szego} is completely integrable
with infinitely many conservation laws
and that analysis of \eqref{Szego} has a strong connection to
the theory of complex variables.

It is shown in \cite{GG} that
\eqref{Szego} is globally
well-posed in $H^\frac{1}{2}_+(\T) :=  \Pi (H^\frac{1}{2}(\T))$
via the energy method (for $ s> \frac{1}{2}$ - the argument for $ s= \frac{1}{2}$
is more intricate) and the conservation of the $H^\frac{1}{2}_+$-norm.
Our interest is to investigate
if there is any nonlinear smoothing
by considering random initial data $u_0$ of the form 
\begin{equation} \label{IVSZ}
u_0(x) = u_0^\omega(x) = \sum_{n \geq 0} \frac{g_n(\omega)}{\sqrt{1+ |n|^{2\al}}} e^{inx}
\in H^{\al - \frac{1}{2}-}_+(\T) \setminus H^{\al - \frac{1}{2}}_+(\T), \ \ \text{a.s.}
\end{equation}

\noi
for $\al \leq 1$.
First, write \eqref{Szego} in the integral formulation:
\begin{equation} \label{Szego2}
u(t) = u_0 - i \mathcal{N}(u, u, u)(t),
\end{equation}

\noi
where $\mathcal{N}(\cdot, \cdot, \cdot)$ is given by
\begin{equation} \label{SzDuhamel}
\mathcal{N}(u_1, u_2, u_3)(t) = 
\int_0^t \Pi(u_1 \cj{u_2} u_3) (t') dt'.
\end{equation}

\noi
In the following, we consider the {\it second iterate}:
\begin{equation} \label{Sz2}
z(t) = u_0 - i \mathcal{N}(u_0, u_0, u_0)(t),
\end{equation}

\noi
where $u_0$ is as in \eqref{IVSZ}.

\begin{proposition} \label{PROP:Szego}
Let $u_0$ be as in \eqref{IVSZ}.

\noi
\textup{(a)} Let $\al \in (\frac{1}{2}, 1]$.
Then for $ s \geq \al - \frac{1}{2}$, 
we have
\[ \big\| \mathcal{N}(u_0, u_0, u_0) \big\|_{C([-T, T];H^s_+)} =\infty, \quad a.s.\]

\noi
for any $T>0$.
In particular, even with $\al = 1$, 
the second iterate $z(t)$ for the cubic Szeg\"o equation \eqref{Szego}
is a.s.~unbounded in $H^\frac{1}{2}_+$. 

\noi \textup{(b)} Let $\al > \frac{1}{2}$.
Then, the second iterate $z(t)$ for the cubic Szeg\"o equation \eqref{Szego}
is a.s.~bounded in $H^s_+$ with $s <\al -\frac{1}{2}-\eps$ for any $\eps > 0$.

\end{proposition}

\noi
Proposition \ref{PROP:Szego} (a)
shows that there is no gain of regularity, even for $\al = 1$.
In view of well-posedness of \eqref{Szego} in $H^\frac{1}{2}_+$, 
we conclude that there is no smoothing upon randomization of initial data.
See also Remark \ref{REM:last}.
This shows that dispersion is closely related to nonlinear smoothing under randomization of initial data.
See Section 5 for details.

\begin{remark} \rm
It is interesting to compare Proposition \ref{PROP:Szego}
with the boundedness of the second iterate for KdV;
in Remark \ref{REM:KDV2}, 
we show that 
the nonlinear part of the second iterate for KdV is bounded in $L^2(\T)$
even for deterministic mean-zero initial data $u_0$ in $ H_0^s(\T)$
as long as $s > -\frac{3}{4}$.

On the one hand, 
the failure of well-posedness 
for the cubic Szeg\"o equation \eqref{Szego} below $H_+^\frac{1}{2}(\T)$ is 
due to the unboundedness of the second iterate below $H_+^\frac{1}{2}(\T)$
(indeed even in $H_+^\frac{1}{2}(\T)$ - see \cite{GG}),
and this problem can not be removed by considering the random initial data of the form \eqref{IVSZ}.
On the other hand, the failure of well-posedness for KdV below $H_0^{-\frac{1}{2}}(\T)$
 is not caused by the second iterate,
 and randomization of initial data with the second iteration argument comes in rescue.
We also point out that 
while the 
 failure of well-posedness of Wick ordered cubic NLS \eqref{NLS0} below $L^2$
is due to the unboundedness of the second iterate below $L^2$,
this problem can be removed by a simple fixed point argument (around the linear solution)
with random initial data.

Lastly, recall an analogous result for the Benjamin-Ono equation by Tzvetkov \cite{TZ3}.
He showed that the second iterate with initial data $u_0$ of the form \eqref{IVKDV}
with $\al = \frac{1}{2}$ is in $H^s(\T) \setminus L^2(\T)$ for any $s < 0$.
(Recall that the threshold regularity of the deterministic well-posedness theory is $L^2(\T)$.
See Molinet \cite{MOLI1, MOLI2}.)
On the one hand, this result shows that there is no nonlinear smoothing under randomization
of initial data.  On the other hand, the second iterate is at least as regular as the random initial data.
Proposition \ref{PROP:Szego} (b) also states that the second iterate is at least as regular as the random initial data.
Although there is no nonlinear smoothing under randomization of initial data, 
these results state that there is still a possibility of constructing solutions
below the deterministic threshold regularities.
However, such a construction is out of reach at this point.
\end{remark}

\medskip

This paper is organized as follows:
We introduce notations in Section 2
and state basic lemmata in Section 3.
In Section 4, we prove Theorem \ref{THM:LWP}
by establishing the nonlinear estimate on the second iteration 
of the integral formulation.
In Section 5, we show that there is no 
extra smoothing for the cubic Szeg\"o equation
upon randomization of initial data.

\section{Notation}
Let $X^{s, b}$ denote the periodic Bourgain space defined in \eqref{Xsb}.
We often use the shorthand notation $\| \cdot\|_{s, b}$
to denote the $X^{s, b}$ norm.
Since the $X^{s, \frac{1}{2}}$ norm fails to control $L^\infty_tH^s_x$ norm,
we introduce a smaller space 
 $Z^{s, b}(\mathbb{T} \times \mathbb{R})$ whose norm is given by 
\begin{equation} \label{Zsb}
\| u\|_{Z^{s, b}(\mathbb{T} \times \mathbb{R})} 
:= \| u\|_{X^{s, b}(\mathbb{T} \times \mathbb{R})}
+ \| u\|_{Y^{s, b-\frac{1}{2}}(\mathbb{T} \times \mathbb{R})}
\end{equation}

\noi
where $\jb{\, \cdot\, } = 1 + |\cdot|$ and 
$\| u\|_{Y^{s, b}(\mathbb{T} \times \mathbb{R})} 
= \|\jb{n}^s \jb{\tau - n^3}^b \ft{u}(n, \tau)\|_{l^2_n L^1_\tau(\mathbb{Z} \times\mathbb{R})}.$
We also define the local-in-time version $Z^{s, b, T }$  on $\T \times [-T, T]$,
by
\[ \|u\|_{Z^{s, b, T }} =  \inf \big\{ \|\wt{u} \|_{Z^{s, b}(\T \times \mathbb{R})}: {\wt{u}|_{[-T, T]} = u}\big\}.\]

\noi
The local-in-time versions of other function spaces are defined analogously.

If a function depends on both $x$ and $t$, we use ${}^{\wedge_x}$ 
(and ${}^{\wedge_t}$) to denote the spatial (and temporal) Fourier transform, respectively.
However, when there is no confusion, we simply use ${}^\wedge$ to denote the spatial Fourier transform,
the temporal Fourier transform, and  the space-time Fourier transform, depending on the context.
For simplicity, we often drop $2\pi$ in dealing with the Fourier transforms.
If a function $f$ is random, we may use the superscript
$f^\omega$ to show the dependence on $\omega$

Lastly,
let $\eta \in C^\infty_c(\mathbb{R})$ be a smooth cutoff function supported on $[-2, 2]$ with $\eta \equiv 1$ 
on $[-1, 1]$
and let $\eta_{_T}(t) =\eta(T^{-1}t)$. 
We use $c,$ $ C$ to denote various constants, usually depending only on $\al$ and $s$.
If a constant depends on other quantities, we will make it explicit.
We use $A\lesssim B$ to denote an estimate of the form $A\leq CB$.
Similarly, we use $A\sim B$ to denote $A\lesssim B$ and $B\lesssim A$
and use $A\ll B$ when there is no general constant $C$ such that $B \leq CA$.
We also use $a+$ (and $a-$) to denote $a + \eps$ (and $a - \eps$), respectively,  
for arbitrarily small $\eps \ll 1$.

\section{Basic lemmata and linear estimates}
We first state several useful lemmata.
See \cite{CO1} for the proofs.
Recall  that by restricting the Bourgain spaces onto a small time interval $[-T, T]$, 
we can gain a small power of $T$ (at a  loss of regularity on $\jb{\tau - n^3} $.)
\begin{lemma} \label{LEM:timedecay}
For $b > b'>0$, we have
\begin{equation} \label{CCtime0}
\|u\|_{X^{s, b', T}} =\|\eta_{_T}u\|_{X^{s, b', T}} 
\lesssim T^{b-b'-} \|u\|_{X^{s, b, T}}.
\end{equation}

\end{lemma}

\noi
The proof basically follows from 
\begin{equation} \label{eq:delta}
\|\ft{\eta_{_{T}}}\|_{L^q_\tau} \sim T^\frac{q-1}{q} \|\ft{\eta}\|_{L^q_\tau} 
\sim T^\frac{q-1}{q},
\end{equation}

\noi
where $\ft{\eta_{_{T}}}(\tau) = T \ft{\eta}(T\tau)$, and interpolation.

\medskip
Next, we present a probabilistic lemma related to the Gaussian
random variables.

\begin{lemma} \label{LEM:prob1}
Let $\eps, \beta > 0$.  Then,  for $T >0 $, we have 
 \begin{equation}
 |g_n(\omega)| \leq C_\eps T^{-\frac{\beta}{2}} \jb{n}^\eps
 \end{equation}

\noi
for all $n \in \mathbb{Z}$ for $\omega$ outside an exceptional set of measure $< e^{-\frac{c}{{T^\beta}}}$.
\end{lemma}


\medskip

Now, we briefly go over the linear  estimates
related to KdV.
Let $S(t) = e^{- t \dx^3  }$ and $T\leq 1$ in the following.
We first present the homogeneous and nonhomogeneous linear estimates.
See \cite{BO1, KPV93} for details.

\begin{lemma} \label{LEM:linear1}
For any $s \in \mathbb{R}$ and $b < \frac{1}{2}$, we have 
$\|  S(t) u_0\|_{X^{s, b, T}} \lesssim T^{\frac{1}{2}-b}\|u_0\|_{H^s}$.
\end{lemma}

\begin{lemma} \label{LEM:linear2}
For any $s \in \mathbb{R}$ and $b \leq \frac{1}{2}$, we have 
\begin{align*} 
 \bigg\|  \int_0^t S(t-t') F(x, t') dt'\bigg\|_{X^{s, b, T}} 
\lesssim \| F \|_{Z^{s, b-1, T}}.
\end{align*}

\noi
Also, we have 
$ \big\|  \int_0^t S(t-t') F(x, t') dt'\big\|_{X^{s, b, T}} 
\lesssim \| F \|_{X^{s, b-1}}$
for $b > \frac{1}{2}$.
\end{lemma}

\noi
The next lemma is the periodic $L^4$-Strichartz estimate for KdV due to Bourgain \cite{BO11}.
\begin{lemma} \label{LEM:L4}
Let $u$ be a function on $\T \times \R$.
Then, we have 
$ \|u\|_{L^4_{x, t}} \lesssim \|u\|_{X^{0, \frac{1}{3}}}.$
\end{lemma}

\section{On the KdV equation} \label{SEC:4}

In this section, we construct local-in-time solutions
to \eqref{KDV} with random initial data of the form \eqref{IVKDV}
for $\al \in ( \al_0, 0]$, where $\al_0$ is as in Theorem \ref{THM:LWP}.

\subsection{Overview; Unboundedness of nonlinear term} \label{SUBSEC:4.1}
By writing  KdV \eqref{KDV} in the Duhamel formulation, we have
\begin{equation} \label{KDVduhamel}
u(t) =  S(t) u_0 + \mathcal{N}(u, u) (t), 
\end{equation} 

\noi
where $S(t) = e^{-t \dx^3}$ and $\mathcal{N}(\cdot, \cdot)$  is  given by
\begin{equation} \label{NN}
\mathcal{N}(u_1, u_2) (t) :=  -\frac{1}{2}\int_0^t S(t - t') \dx (u_1 u_2)(t') d t'.
\end{equation}

It follows from the conservation of mean 
that $u(t)$ has the spatial mean 0 for each $t\in \R$
since $u_0$ has the mean 0.
We use $(n, \tau)$, $(n_1, \tau_1)$, and $(n_2, \tau_2)$ to denote the Fourier variables
for $uu$, the first factor, and the second factor $u$ of $uu$ in $\mathcal{N}(u, u)$, 
respectively.
i.e. we have $n = n_1 + n_2$ and $\tau = \tau_1 + \tau_2$.
By the mean zero assumption on $u$ and 
by the fact that we have $\dx (uu)$ in the definition of $\mathcal{N}(u,u)$, 
we may assume  $n, n_1, n_2 \ne 0$.
We also use the following notation:
\[\s_0 := \jb{\tau - n^3} \text{ and } \s_j := \jb{\tau_j - n_j^3}.\]

\noi
One of the main ingredients is the observation due to Bourgain \cite{BO11}:
\begin{equation} \label{Walgebra}
n^3 - n_1^3 - n_2^3 = 3 n n_1 n_2, \ \text{for } n = n_1 + n_2,
\end{equation}

\noi
which in turn implies that 
\begin{align} \label{MAXMAX}
\MAX& := \max( \s_0, \s_1, \s_2) \notag \\
& \gtrsim \jb{(\tau-n^3) - (\tau_1 - n_1^3) - (\tau_2- n_2^3)} \sim \jb{n n_1 n_2}.
\end{align}

\noi
 This estimate \eqref{MAXMAX} played a crucial role 
in establishing the bilinear estimate \eqref{KPVbilinear}.

\medskip

If we were to proceed as in \cite{BO7, CO1}, 
then we would need to estimate $\|\mathcal{N}(u_1, u_2)\|_{
Z^{-\frac{1}{2}, \frac{1}{2}, T}}$,
assuming that $u_j$ is either of type
\begin{itemize}
\item[(I)] random, less regular:
\begin{equation}  \label{P type 1}
u_j(t) =  \eta_{_T} (t)\sum_{n\ne 0 } g_n(\omega) e^{i(nx + n^3 t)}, \ \ \text{or}
\end{equation}

\item[(II)] deterministic, smoother:
$u_j = v \text{ with }\|v\|_{Z^{-\frac{1}{2}, \frac{1}{2}, T}} \leq 1.$
\end{itemize}

\noi
In the following, we 
show that  estimate on 
$\|\mathcal{N}(u_1, u_2)\|_{Z^{-\frac{1}{2}, \frac{1}{2}, T}}$
fails to hold
by considering the case when both $u_1$ and $u_2$ are of type (I).

By computing the Duhamel term explicitly (see \eqref{Duhamel} below),
one of the main contributions to 
$\|\mathcal{N}(u_1, u_2)\|_{Z^{-\frac{1}{2}, \frac{1}{2}, T}}$
is given by  
$\| \dx(u_1u_2)\|_{X^{-\frac{1}{2}, -\frac{1}{2}, T}}$.
Now, assume that $u_1$ and $u_2$ are of type (I).
i.e., we have $\ft{u}_j (n_j, \tau_j) = \ft{\eta}_{_T}(\tau_j - n_j^3) g_{n_j}$, $j = 1, 2$.

For simplicity of the presentation, 
we remove the time cutoff $\eta$ in \eqref{P type 1}.
Then, we have $\ft{u}_j (n_j, \tau_j) = \dl(\tau_j - n_j^3) g_{n_j}$, $j = 1, 2$.
Hence, from \eqref{MAXMAX} with $\s_1, \, \s_2 \sim1$, 
we have $\s_0 = \MAX \sim\jb{n n_1 n_2}$.
Then, we have
\begin{align} \label{P nonlin}
\notag
\| \dx(u_1u_2) & \|_{X^{-\frac{1}{2}, -\frac{1}{2}}}^2
 = \sum_n \int  \bigg|\frac{|n|\jb{n}^{-\frac{1}{2}}}{\jb{\tau-n^3}^\frac{1}{2}}
\sum_{n = n_1 + n_2} \intt_{\tau=\tau_1 + \tau_2}
\ft{u}_1 (n_1, \tau_1) \ft{u}_2 (n_2, \tau_2)d\tau_1\bigg|^2
d\tau\\
& \notag
\sim \sum_n \int  \bigg|\sum_{n = n_1 + n_2} \int_{\tau=\tau_1 + \tau_2}
\frac{\dl(\tau_1 - n_1^3)g_{n_1}}{|n_1|^\frac{1}{2}}
\frac{\dl(\tau_2 - n_2^3)g_{n_2}}{|n_2|^\frac{1}{2}}
d\tau_1\bigg|^2
d\tau \\
& = \sum_n \int  
\sum_{\substack{n = n_1 + n_2\\\hphantom{llll}=m_1+m_2}} 
\intt_{\tau=\tau_1 + \tau_2}
\prod_{j = 1}^2 \frac{\dl(\tau_j - n_j^3)g_{n_j}}{|n_j|^\frac{1}{2}}
d\tau_1 
\intt_{\tau=\wt{\tau}_1 + \wt{\tau}_2}
\prod_{k = 1}^2 \frac{\dl(\wt{\tau}_k - m_k^3)\cj{g_{m_k}}}{|m_k|^\frac{1}{2}}
d\wt{\tau}_1.
\end{align}

\noi
We have  nontrivial contribution in \eqref{P nonlin}
only when $n_1^3 + n_2^3 = m_1^3 + m_2^3$
and $n_1 + n_2 = m_1 + m_2$.
Hence, we have $\{n_1, n_2\} = \{m_1, m_2\}$.
Therefore, we have
\begin{align*}
\eqref{P nonlin}
\sim \sum_{n} \sum_{n = n_1 + n_2} \frac{|g_{n_1}|^2}{|n_1|}\frac{|g_{n_2}|^2}{|n_2|}
=  \sum_{n_1} \frac{|g_{n_1}|^2}{|n_1|} \sum_{n_2}\frac{|g_{n_2}|^2}{|n_2|}
= \infty, \ \text{ a.s.}
\end{align*}

\noi
The last equality holds from the following.
Let $F_j (\omega) := 2^{-j}\sum_{|n|\sim 2^j} |g_n(\omega)|^2$.
Then, $F_j$ converges to $\text{Var}(g_n) = 2$ a.s.
by strong law of large numbers.
Hence, the tails of the above sums do not converge to 0.

The above computation involving the Dirac delta function
is somewhat formal.
It can be made rigorous
by using a smooth cutoff function $\eta_{_T}$.
However, we omit details.
As a conclusion, we see that 
a simple application of the ideas from \cite{BO7, CO1}
fails for KdV.
In the remaining part of this section,
we construct local-in-time solutions
by adapting the second iteration argument \cite{BO3, OH6}
in the probabilistic setting.

\begin{remark} \label{REM:KDV2} \rm
Consider the second iterate for the Duhamel formulation \eqref{KDVduhamel}
of KdV:
\begin{equation} \label{KDV2}
 z(t) =  S(t) u_0 + \mathcal{N}(S(t)u_0, S(t)u_0) (t), 
\end{equation}

\noi
where $\mathcal{N}(\cdot, \cdot)$ is defined in \eqref{NN}
and $u_0 \in H^s_0(\T)$.
In the following, we show that the nonlinear part of the second iterate
is bounded even in $L^2(\T)$ as long as 
$u_0$ is in  $H^s_0(\T)$ for $s> -\frac{3}{4}$.
In particular, this shows that 
the failure of (analytic) well-posedness for KdV below $H_0^{-\frac{1}{2}}(\T)$
 is not due to  the second iterate.

Fix $t>0$. Then, we have
\begin{align*}
\big\|\mathcal{N}(S(t)u_0, S(t)u_0) (t) \big\|_{L^2}
\sim \bigg(\sum_{n} 
\bigg|\int_0^t e^{-it'n^3} in \sum_{n = n_1 + n_2} \prod_{j = 1}^2 e^{it'n_j^3} \ft{u}_0(n_j) dt'\bigg|^2 \bigg)^\frac{1}{2}.
\end{align*}

\noi Since $u_0$ has mean zero on $\T$, we assume $n_1, n_2 \ne 0$ in the following.
Moreover, we assume $n \ne 0$ thanks to the derivative in the nonlinearity $\dx(u_1 u_2)$.
First integrate in $t'$ with \eqref{Walgebra}.
Then, by Young's inequality followed by H\"older's inequality, we have
\begin{align*}
\big\|\mathcal{N} (S(t)u_0, & S(t)u_0)  (t) \big\|_{L^2}
 \lesssim \bigg(\sum_{n} 
\bigg| in \sum_{n = n_1 + n_2} \prod_{j = 1}^2  \ft{u}_0(n_j)  \int_0^t e^{-it'(n^3-n_1^3 - n_2^3)} dt' \bigg|^2 \bigg)^\frac{1}{2}\\
& \lesssim \bigg(\sum_{n} 
\Big|  \sum_{n = n_1 + n_2} \prod_{j = 1}^2  |n_j|^{-1} |\ft{u}_0(n_j)|   \Big|^2 \bigg)^\frac{1}{2}\\
&  \lesssim \prod_{j = 1}^2 \big\| \jb{n_j}^{-1} |\ft{u}_0(n_j)| \big\|_{l^\frac{4}{3}_{n_j}} 
 \leq \prod_{j = 1}^2 \| \jb{n_j}^{-1-s} \|_{l^{4}_{n_j}}  \|\jb{n_j}^{s} \ft{u}_0(n_j)\|_{l^2_{n_j}}
\lesssim \|u_0\|_{H^s}^2
\end{align*}

\noi
as long as $4(-1-s) <-1$, i.e. $s>-\frac{3}{4}$.
By considering the random initial data $u_0$ of the form \eqref{IVKDV}, 
we can also show that the nonlinear part of the second iterate
is a.s.~bounded in $L^2$ as long as $\al > -\frac{1}{2}$, namely $s> -1$.
We omit the details.

\end{remark}

\subsection{Nonlinear analysis via second iteration}
\label{SUBSEC:4.2}
First, we briefly go over Bourgain's argument in \cite{BO3}.
Define 
\begin{equation} \label{AJJ}
A_j = \{(n, n_1, n_2, \tau, \tau_1, \tau_2) \in \mathbb{Z}^3 \times \R^3:
\s_j = \MAX\},
\end{equation}

\noi
and let $\mathcal{N}_j(u, u)$ denote the contribution of $\mathcal{N}(u, u)$ on $A_j$.
Then, \eqref{KDVduhamel} can be written as 
\begin{equation} \label{KDVduhamel2}
u(t) =  S(t) u_0 + \mathcal{N}_0(u, u) (t)+
\mathcal{N}_1(u, u) (t) + \mathcal{N}_2(u, u) (t).
\end{equation}

\noi
By the standard bilinear estimate with Lemma \ref{LEM:L4} as in \cite{BO11}, \cite{KPV4}, 
we have
\begin{align} \label{N_0}
\|\mathcal{N}_0(u, u)\|_{{-\frac{1}{2} + \dl, \frac{1}{2}-\dl}}
\leq o(1)\|u\|^2_{{-\frac{1}{2} - \dl, \frac{1}{2}-\dl}},
\end{align}

\noi
where $o(1) = T^\theta$ with $\theta > 0$ by considering the estimate 
on a short time interval $[-T, T]$.
See (2.17), (2.26), and (2.68) in \cite{BO3}.
Here, we abuse the notation and use $\|\cdot\|_{s, b} = \|\cdot\|_{X^{s, b}}$ 
to denote the local-in-time version as well.
Note that the temporal regularity \[b = \tfrac{1}{2} - \dl < \tfrac{1}{2}.\]

\noi
This allowed us to gain the spatial regularity by $2\dl$ in \eqref{N_0}.
Clearly, we can not expect to do the same for $\mathcal{N}_1(u, u)$.
(By symmetry, we do not consider $\mathcal{N}_2(u, u)$ in the following.)
When $b < \frac{1}{2}$, the bilinear estimate \eqref{KPVbilinear} is known to fail for any $s \in \R$
 due to the contribution from $\mathcal{N}_1(u, u)$. 
See \cite{KPV4}.
Following the notation in \cite{BO3}, 
let 
\begin{equation}\label{Isb}
I_{s, b} = \|\mathcal{N}_1(u, u) \|_{X^{s, b}} \
\text{ and }\ a := \tfrac{1}{2} -\dl < \tfrac{1}{2}.
\end{equation}

\medskip
\noi
{\bf Main goal:} Estimate the Duhamel term $\mathcal{N}_1 (u, u)$
in $X^{s, b}$ with
\begin{equation} \label{s1} s :=  -a = -\tfrac{1}{2} + \dl >-\tfrac{1}{2}, \quad \text{and} \quad 
b :=  1-a = \tfrac{1}{2}+\dl >\tfrac{1}{2}, 
\end{equation}

\noi
assuming that $u$ is bounded only in $X^{\wt{s}, \wt{b}}$
with 
\begin{equation} \label{s2} \wt{s} :=  -(1-a) = -\tfrac{1}{2} - \dl <-\tfrac{1}{2} , \quad \text{and} \quad 
\wt{b} := a = \tfrac{1}{2}-\dl <\tfrac{1}{2}. 
\end{equation}

\medskip

By Lemma \ref{LEM:linear2} and duality 
with $\|d(n, \tau)\|_{L^2_{n, \tau}} \leq 1$, we have 
\begin{align} \label{eq:I1}
I_{-a, 1-a} & = \|\mathcal{N}_1(u, u) \|_{-a, 1-a}\\
& \lesssim \sum_{\substack{n, n_1\\n = n_1 + n_2}} \intt_{\tau = \tau_1 + \tau_2} d\tau d\tau_1
\frac{\jb{n}^{1-a}d(n, \tau)}{\s_0^{a }} \ft{u}(n_1, \tau_1) \frac{\jb{n_2}^{1-a}c(n_2, \tau_2)}{\s_2^a},
\notag
\end{align}

\noi
where  
\begin{equation} \label{CN}
c(n_2, \tau_2) = \jb{n_2}^{-(1-a)}\s_2^{a}\, \ft{u}(n_2, \tau_2) 
\text{ so that }
\|c\|_{L^2_{n, \tau}} = \|u\|_{-(1-a), a} = \|u\|_{-\frac{1}{2}-\dl, \frac{1}{2}-\dl}.
\end{equation}

\noi
The main idea here is to consider the second iteration, 
i.e. substitute \eqref{KDVduhamel} for  $\ft{u}(n_1, \tau_1)$ in \eqref{eq:I1}, 
thus leading to a trilinear expression, i.e.
write 
$\mathcal{N}_1(u, u)$ as 
\begin{align} \label{KDVduhamel3}
\mathcal{N}_1(u, u) = 
\mathcal{N}_1(S(t) & u_0, u)
 + \mathcal{N}_1(\mathcal{N}(u, u), u).
\end{align}

\noi
Applying the second iteration on the second argument of $\mathcal{N}_2(u, u)$, 
we can write \eqref{KDVduhamel} and \eqref{KDVduhamel2} as
\begin{align} \label{KDVduhamel4} \notag
u(t) =  S(t) u_0 & + \mathcal{N}_0(u, u) (t)
 + \mathcal{N}_1(S(t)  u_0, u)(t)\\
& + \mathcal{N}_1(\mathcal{N}(u, u), u) (t) 
+ \mathcal{N}_2(u, S(t)  u_0)(t)
+ \mathcal{N}_2(u, \mathcal{N}(u, u)) (t).
\end{align}

\noi
Since $\s_1 = \MAX \gtrsim \jb{n n_1n_2}\gg1$ on $A_1$, 
we can assume that 
\begin{equation} \label{eq:u_1}
\ft{u}(n_1, \tau_1) = \big(\mathcal{N}(u, u)\big)^\wedge(n_1, \tau_1)
\sim \frac{|n_1|}{\s_1} \sum_{n_1 = n_3 + n_4} \intt_{\tau_1 = \tau_3 + \tau_4}
\ft{u}(n_3, \tau_3)\ft{u}(n_4, \tau_4) d\tau_4.
\end{equation}

\noi
Namely, we can assume that the contribution to $\ft{u}(n_1, \tau_1)$ 
from the linear part $S(t) u_0$  of \eqref{KDVduhamel2} is negligible
since we have $\s_1 \sim 1$ for the linear part.
Moreover, by the standard computation \cite{BO11}, we have
\begin{align} \label{Duhamel}
 \mathcal{N}(u, u)(x, t) 
& = -i \sum_{k= 1}^\infty \frac{i^k t^k}{k!}  
\sum_{n \ne 0} e^{i(nx + n^3t)} \int \eta(\ld - n^3) (\ld - n^3)^{k-1}\ft{\dx u^2}(n, \ld) d\ld \notag \\
& \hphantom{X}+ i \sum_{n \ne 0} e^{inx} \int 
\frac{\big(1-\eta\big)(\tau - n^3)}{\tau - n^3}  
\ft{\dx u^2}(n, \tau) e^{i \tau t } d \tau \notag \\
& \hphantom{X}+ i  \sum_{n \ne 0} e^{i(nx + n^3 t)} \int 
\frac{\big(1-\eta\big)(\ld - n^3)}{\ld - n^3}  \ft{\dx u^2}(n, \ld) d \ld \notag \\
& =:  \mathcal{M}_1(u, u)(x, t) + \mathcal{M}_2(u, u)(x, t)+ \mathcal{M}_3(u, u)(x, t).
\end{align}

\noi
Note that 
$(\mathcal{M}_1(u, u))^\wedge(n_1, \tau_1)$
and $(\mathcal{M}_3(u, u))^\wedge(n_1, \tau_1)$
are distributions supported on $\{\tau_1 -n_1^3 = 0\}$.
i.e. $\s_1 \sim 1$.
Hence, the only contribution for the second iteration on $A_1$ comes from 
$\mathcal{M}_2(u, u)$
whose Fourier transform is given in 
\eqref{eq:u_1}.
This shows the validity of the assumption \eqref{eq:u_1}.

The $\s_1 $ appearing in the denominator
allows us to cancel $\jb{n}^{1-a}$ and $\jb{n_2}^{1-a}$ in the numerator in \eqref{eq:I1}.
Then, $I_{-a, 1-a}$ can be estimated by 
\begin{align} \label{eq:I2}
\lesssim \sum_{\substack{n = n_1 + n_2\\n_1 = n_3 + n_4}} 
\intt_{\substack{\tau = \tau_1 + \tau_2\\\tau_1 = \tau_3 + \tau_4}} 
\frac{\jb{n}^{1-a}d(n, \tau)}{\s_0^{a }} \frac{|n_1|}{\s_1} \, \ft{u}(n_3, \tau_3)\ft{u}(n_4, \tau_4)
\frac{\jb{n_2}^{1-a}c(n_2, \tau_2)}{\s_2^a}
.
\end{align}

\noi
The argument was then divided into several cases, 
depending on the sizes of $\s_0, \cdots, \s_4$.
Here, the key algebraic relation is
\begin{equation}\label{algebra2}
n^3 - n_2^3 - n_3^3 - n_4^3 = 3(n_2+ n_3)(n_3+ n_4)(n_4+ n_2), \  \text{ with } n = n_2 + n_3 + n_4.
\end{equation}

\noi
Then, Bourgain proved -see (2.69) in \cite{BO3}-
\begin{equation} \label{N_1}
I_{-a, 1-a} \leq o(1)\|u\|_{-(1-a), a} I_{-a, 1-a}  
+ o(1) \|u\|^3_{-(1-a), a} + o(1) \|u\|_{-(1-a), a},
\end{equation}

\noi
{\it assuming} the a priori estimate \eqref{BOO}: $|\ft{u}(n, t)| <C$ for all $n\in \mathbb{Z}$, $t \in\R$.
Indeed, the estimates involving the first two terms on the right hand side of  \eqref{N_1}
were obtained without \eqref{BOO}, and
{\it only} the last term in \eqref{N_1} required \eqref{BOO}, 
-see ``Estimation of (2.62)'' in \cite{BO3}-, 
which was then used to deduce
\begin{equation} \label{eq:apriori}
\|\ft{u}(n, \cdot)\|_{L^2_\tau} < C.
\end{equation}

\noi
The a priori estimate \eqref{BOO} is derived via
the isospectral property of the KdV flow
and is false for a general function in $X^{-(1-a), a}$.
(It is here that the smallness of the total variation $\|\mu\|$ is used in \cite{BO3}.)

\medskip

Our goal is to carry out a similar analysis on the second iteration 
{\it without} the a priori estimates \eqref{BOO} and \eqref{eq:apriori} coming from the complete integrability of KdV.
We achieve this goal by exhibiting 
nonlinear smoothing under randomization on initial data.
In the following, we take 
\begin{equation} \label{s3}\al > -\dl = a - \tfrac{1}{2},
\end{equation}

\noi
where $\dl > 0$ is as in \eqref{Isb}, \eqref{s1}, and \eqref{s2}.
Then,  the initial data
\begin{equation}\label{Q IV}
u_0(x) = u_0^\omega(x) = \sum_{n\ne 0} \frac{g_n(\omega)}{|n|^\al} e^{inx}
\end{equation}

\noi
belongs to $H^{\wt{s}}(\T) = H^{-\frac{1}{2}-\dl}(\T)$, a.s.

Note that we can use the estimates on 
 $\mathcal{N}_1(\mathcal{N}(u, u), u)$ from \cite{BO3} 
{\it except} when the a priori bound \eqref{BOO} was assumed.
i.e. we need to estimate the contribution from (2.62) in \cite{BO3}:
\begin{equation}\label{262}
R_a (u_2, u_3, u_4):= 
\sum_{n} \intt_{\tau = \tau_2 +  \tau_3 + \tau_4} 
\chi_B \frac{d(n, \tau)}{\jb{n}^{1+a}\s_0^{a }} 
 \ft{u}_2(-n, \tau_2) \ft{u}_3(n, \tau_3)\ft{u}_4(n, \tau_4)d\tau_2 d\tau_3d\tau_4,
\end{equation}

\noi
where $\|d(n, \tau)\|_{L^2_{n, \tau}} \leq 1$
and $B = \{ \s_0, \s_2, \s_3, \s_4 < |n|^\g\}$ with some small parameter $\g>0$. 
This corresponds to the case $n_2 = -n$ and $n_3 = n_4 = n$
in \eqref{eq:I2} after some reduction.
In our analysis, we directly estimate $R_a(u_2, u_3, u_4)$,
assuming that $u_j$ is either of type 
\begin{itemize}
\item[(I)] linear part:  random, less regular 
\begin{equation*} 
u_j(t) = \eta_{_T}(t) \sum_{n\ne 0 } \frac{g_n(\omega)}{|n|^\al} e^{i(nx + n^3 t)}, \ \ \text{or}
\end{equation*}

\item[(II)]  nonlinear part: deterministic, and (expected to be) smoother
\[u_j = \mathcal{N}(u): = \mathcal{N}(u, u)\ \text{ (to be bounded in  $X^{-\frac{1}{2}+\dl, \frac{1}{2}-\dl, T}$)}.\] 
\end{itemize}

\noi
In \cite{BO3}, this parameter $\g = \g(a)$,  subject to the conditions (2.43) and (2.60) in \cite{BO3},  
played a certain role in estimating $R_a$
along with the a priori bound \eqref{BOO}.
However, it plays no role in our analysis.
Before proceeding further, we record the conditions on the values of $a$ and $\g$ from \cite{BO3}:\footnote{Note that
$\al$ in \cite{BO3} is $a$ in this paper.}
\begin{align} \label{aa}
 a > \tfrac{7}{18}, \quad 
 \g > \frac{2(1-2a)}{a-1/3}, \quad \text{and} \quad
 2(1-a) + \g < 2.
\end{align}

\noi
From the last two conditions, we obtain 
the quadratic inequality $a^2 + \tfrac{5}{3} a - 1> 0$.
By choosing $a> a_0$, where 
$a_0 = -\tfrac{5}{6} +\tfrac{\sqrt{61}}{6} \thickapprox 0.4684$
is as in \eqref{a0} in the statement of Theorem \ref{THM:LWP}, 
we guarantee that all the three conditions in \eqref{aa}
are satisfied.
Hence, we assume that $a \in (a_0, \frac{1}{2})$ in the following.
In particular, this implies that 
\begin{equation} \label{al}
\al  >\al_0 :=  a_0 - \tfrac{1}{2} \thickapprox -0.0316
\end{equation}

\noi
and $\dl < -\al_0 \thickapprox 0.0316$
such that \eqref{s3} holds.

\medskip
\noi
$\bullet$ {\bf Case 1:} all type (II).
\quad By Cauchy-Schwarz and Young's inequalities, we have
\begin{align}
\eqref{262} & \leq \sum_n \| d(n, \cdot)\|_{L^2_\tau}
\jb{n}^{-1-a} \| \ft{\mathcal{N}(u)}(-n, \tau_2)\|_{L^{\frac{6}{5}}_{\tau_2}}
\| \ft{\mathcal{N}(u)}(n, \tau_3)\|_{L^{\frac{6}{5}}_{\tau_3}}
\| \ft{\mathcal{N}(u)}(n, \tau_4)\|_{L^{\frac{6}{5}}_{\tau_4}} \notag
\intertext{By H\"older inequality (with appropriate $\pm$ signs) and the fact that $-1-a \leq -3a$,}
 & \leq \sum_n \| d(n, \cdot)\|_{L^2_\tau}
\prod_{j = 2}^4 \jb{n}^{-a}\|\s_j^{-\frac{1}{3}-}\|_{L^3_{\tau_j}}
\| \s_j^{\frac{1}{3}+} \ft{\mathcal{N}(u)}(\pm n, \tau_j)\|_{L^{2}_{\tau_j}} \notag \\
&  \leq T^{\frac{1}{2}-3\dl -}\| d(\cdot, \cdot)\|_{L^2_{n, \tau}}
\|\mathcal{N}(u)\|_{X^{-a, a}_{6, 2}}^3
\leq T^{\frac{1}{2}-3\dl -} \|\mathcal{N}(u)\|_{X^{-\frac{1}{2}+\dl, \frac{1}{2}-\dl}}^3, \label{RA}
\end{align}

\noi
where the last two inequalities follows 
from Lemma \ref{LEM:timedecay}
by choosing  $a > \frac{1}{3}$.

\medskip

In the following, fix small $\eps > 0$ and $\beta > 0$.

\noi
$\bullet$ {\bf Case 2:} all type (I).
\quad
By Lemma \ref{LEM:prob1}, we have
$|g_n(\omega)| \leq  C T^{-\frac{\beta}{2}} \jb{n}^\eps$
 outside an exceptional set of measure $<e^{-\frac{c}{T^\beta}}$.
Then, by Cauchy-Schwarz, Young's inequalities,
 and Lemma \ref{LEM:timedecay} with \eqref{s3}, we have
\begin{align}
\eqref{262} 
& \lesssim T^{-\frac{3\beta}{2}} \sum_{n} 
\int d(n, \tau)\jb{n}^{-\frac{3}{2}+\dl + 3\eps - 3\al} \notag \\
 & \hphantom{XXXXX}\times
 \bigg(\intt_{\tau = \tau_2 +  \tau_3 + \tau_4} 
 \ft{\eta}_{_T}(\tau_2+n^3) \ft{\eta}_{_T}(\tau_3-n^3)\ft{\eta}_{_T}(\tau_4-n^3)d\tau_2 d\tau_3\bigg) d\tau \notag \\
& \leq T^{-\frac{3\beta}{2}}\sum_n \| d(n, \cdot)\|_{L^2_\tau}
\jb{n}^{-\frac{3}{2}+4\dl + 3\eps} \|\ft{\eta}_{_T}*\ft{\eta}_{_T}*\ft{\eta}_{_T}\|_{L^{2}_{\tau}}
\lesssim T^{\frac{1}{2}-\frac{3\beta}{2}-}, \label{RA2}
\end{align}

\noi
as long as $\dl < \frac{1}{4}-\frac{3}{4}\eps$.

\medskip
\noi
$\bullet$ {\bf Case 3:} two type (I), one type (II).
\quad
Without loss of generality, assume that $u_2, u_3$ are of type (I),
and that $u_4$ is of type (II).
By Lemma \ref{LEM:prob1}, 
Cauchy-Schwarz, Young, H\"older inequalities,
 and Lemma \ref{LEM:timedecay}, we have
\begin{align}
\eqref{262} &  \leq T^{-\beta} \sum_n \| d(n, \cdot)\|_{L^2_\tau}
\jb{n}^{-\frac{3}{2}+\dl +2\eps-2\al} \notag \\
& \hphantom{XXXXXX}\times
 \bigg\|\intt_{\tau = \tau_2 +  \tau_3 + \tau_4} 
\ft{\eta}_{_T}(\tau_2+n^3) \ft{\eta}_{_T}(\tau_3-n^3)\ft{\mathcal{N}(u)}(n, \tau_4)d\tau_2 d\tau_3\bigg\|_{L^2_\tau}\notag\\
& \leq T^{-\beta} \|d\|_{L^2_{n, \tau}}
\|\ft{\eta}_{_T}\|^2_{L^\frac{6}{5}}
\Big\| \jb{n}^{-\frac{3}{2}+3\dl +2\eps}  \|\ft{\mathcal{N}(u)}(n, \tau)\|_{L^\frac{6}{5}_\tau}\Big\|_{L^2_n}\notag\\
& \lesssim T^{\frac{1}{3}-\beta-}
\Big\| \jb{n}^{-\frac{3}{2}+3\dl +2\eps}  \|\s_4^{-\frac{1}{3}-}\|_{L^3_\tau}
\|\s_4^{\frac{1}{3}+} \ft{\mathcal{N}(u)}(n, \tau)\|_{L^2_\tau}\Big\|_{L^2_n}\notag \\
& \lesssim T^{\frac{1}{2}-\dl-\beta-}
\|\mathcal{N}(u)\|_{X^{-\frac{1}{2}+\dl,\frac{1}{2}-\dl}}  \label{RA3}
\end{align}

\noi
outside an exceptional set of measure $<e^{-\frac{c}{T^\beta}}$
as long as $\dl \leq \frac{1}{2}-\eps$.

\medskip
\noi
$\bullet$ {\bf Case 4:} one type (I), two type (II).
\quad
Without loss of generality, assume that $u_2$ is of type (I),
and that $u_3, u_4$ are of type (II).
By Lemma \ref{LEM:prob1}, 
Cauchy-Schwarz, Young, H\"older inequalities,
 and Lemma \ref{LEM:timedecay}, we have
\begin{align}
\eqref{262} & \leq  T^{-\frac{\beta}{2}} \sum_n \| d(n, \cdot)\|_{L^2_\tau}
\jb{n}^{-\frac{3}{2}+\dl +\eps-\al} \notag \\
& \hphantom{XXXXXX}\times
 \bigg\|\intt_{\tau = \tau_2 +  \tau_3 + \tau_4} 
\ft{\eta}_{_T}(\tau_2+n^3) \ft{\mathcal{N}(u)}(n, \tau_3)\ft{\mathcal{N}(u)}(n, \tau_4)d\tau_2 d\tau_3\bigg\|_{L^2_\tau}\notag \\
& \leq T^{-\frac{\beta}{2}} \|d\|_{L^2_{n, \tau}}
\|\ft{\eta}_{_T}\|_{L^\frac{6}{5}}
\Big\| \jb{n}^{-\frac{3}{2}+2\dl +\eps} 
\prod_{j = 3}^4
\|\s_j^{-\frac{1}{3}-}\|_{L^3_\tau}
\|\s_j^{\frac{1}{3}+} \ft{\mathcal{N}(u)}(n, \tau_j)\|_{L^2_{\tau_j}}
\Big\|_{L^2_n} \notag \\
& \lesssim T^{\frac{1}{6}-\frac{\beta}{2}-}
\|\mathcal{N}(u)\|_{X^{-\frac{1}{2}+\dl,\frac{1}{3}+}}^2 
\lesssim T^{\frac{1}{2}-2\dl-\frac{\beta}{2}-}\|\mathcal{N}(u)\|_{X^{-\frac{1}{2}+\dl,\frac{1}{2}-\dl}}^2 \label{RA4}
\end{align}

\noi
outside an exceptional set of measure $<e^{-\frac{c}{T^\beta}}$
as long as $\eps \leq \frac{1}{2}$.

\medskip
Lastly, we point out that 
the estimates hold even
if we restrict the initial data to be supported on $\{ |n| \leq N\}$,
independent of $N$.
Moreover, we can gain extra power of $N^{0-}$
in Cases 2--4, 
if we restrict the initial data to be supported on $\{ |n| > N\}$.

\subsection{Local well-posedness}
\label{SUBSEC:4.3}
Consider initial data $u_0^N$ of the form
\begin{equation}\label{Q IVN}
u_0^N(x) = \sum_{1 \leq |n|\leq N} \frac{g_n(\omega)}{|n|^\al} e^{inx}.
\end{equation}

\noi
Then, for each $N$, there exists
a unique global solution $u^N \in C(\R; H^{s}(\T))$ for any $s \geq -\frac{1}{2}$ a.s. in $\omega$.
Define $\G^N = \G^N_{u_0^N}$
by
\[ \G^N v = \G^N_{u_0^N} v := S(t) u_0^N +\mathcal{N}(v, v).\]

\noi
Then, $u^N = \G^N u^N$. 
We set $\mathcal{N}^N := \mathcal{N}(u^N) = \mathcal{N}(u^N, u^N)$.

Now, we put all the {\it a priori} estimates together.
Note that all the implicit constants are independent of $N$.
Also, when there is no superscript $N$, 
it means that $N = \infty$.
In the following, $C_j$, $\theta_j$, and $\eps_j$ denote
positive constants.

Fix $a = \frac{1}{2} - \dl >a_0$ as in \eqref{Isb},
where $a_0$ is defined in \eqref{a0}.
From Lemma \ref{LEM:linear1}, we have 
\begin{equation} \label{WT0}
\| S(t) u_0^N\|_{X^{s, b, T}} \leq C_1 \|u_0^N\|_{H^s}
\end{equation}

\noi
for any $s, b\in \R$ with $C_1 =  C_1(b)$.
In particular, by taking $b >\frac{1}{2}$,
we see that $S(t)u_0$ is continuous on $[-T, T]$
with values in $H^s$.
From the definition of $\mathcal{N}_j(\cdot, \cdot)$, \eqref{N_0}, and \eqref{Isb},
we have
\begin{align}  \label{WT1}
\|\mathcal{N}(u^N, u^N)\|_{X^{-a, a, T}}
\leq C_2 T^{\theta_1} \|u^N\|^2_{X^{-(1-a), a, T}}
+ 2 I^N_{-a, a}.
\end{align}

\noi
From \eqref{N_1} and \eqref{RA}--\eqref{RA4}, 
there exists a set $\Omega^{(1)}_T$ with $P\big((\Omega^{(1)}_T)^c\big) < e^{-\frac{c}{T^\beta}}$
such that we have
\begin{align*}
I^N_{-a, 1-a} \leq C_3\big( T^{\theta_2} 
& \|u^N\|_{X^{-(1-a), a, T}}  I^N_{-a, 1-a}  
+ T^{\theta_3} \|u^N\|_{X^{-(1-a), a, T}}^3 
+ T^{\theta_4} \|\mathcal{N}^N\|_{X^{-a, a, T}}^3 \\
& +T^{\theta_5} \|\mathcal{N}^N\|_{X^{-a, a, T}}^2 
+T^{\theta_6} \|\mathcal{N}^N\|_{X^{-a, a, T}} 
+T^{\theta_7} 
\big)
\end{align*}

\noi
on $\Omega^{(1)}_T$.
Note that the choice of $\Omega^{(1)}_T$ is independent of $N$.
For fixed $R > 0$, 
choose $T>0$ sufficiently small such that $C_3T^{\theta_2} R \leq \frac{1}{2}$.
Then, we have 
\begin{align} \label{WT4}
I^N_{-a, 1-a} \leq 2C_3\big( 
&  T^{\theta_3} \|u^N\|_{X^{-(1-a), a, T}}^3 
+ T^{\theta_4} \|\mathcal{N}^N\|_{X^{-a, a, T}}^3 \notag \\
& +T^{\theta_5} \|\mathcal{N}^N\|_{X^{-a, a, T}}^2 
+T^{\theta_6} \|\mathcal{N}^N\|_{X^{-a, a, T}} 
+T^{\theta_7} 
\big)
\end{align}

\noi
for $\|u^N\|_{X^{-(1-a), a, T}} \leq R$.
From \eqref{WT0}--\eqref{WT4}, we have
\begin{align} \label{GAMMA1}
\| u^N \|_{X^{-(1-a), a, T}}
   \leq  C_1 & \|u_0^N\|_{H^{-(1-a)}}
 + \tfrac{1}{2} C_2 T^{\theta_1} \|u^N\|^2_{X^{-(1-a), a, T}} \notag \\
& + 2C_3\big( 
  T^{\theta_3} \|u^N\|_{X^{-(1-a), a, T}}^3 
+ T^{\theta_4} \|\mathcal{N}^N\|_{X^{-a, a, T}}^3  \\
& +T^{\theta_5} \|\mathcal{N}^N\|_{X^{-a, a, T}}^2 
+T^{\theta_6} \|\mathcal{N}^N\|_{X^{-a, a, T}} 
+T^{\theta_7} \big) \notag
\end{align}

\noi
and
\begin{align} \label{GAMMA5}
\| \mathcal{N}^N \|_{X^{-a, a, T}}
   \leq  & \ 
 \tfrac{1}{2} C_2 T^{\theta_1} \|u^N\|^2_{X^{-(1-a), a, T}} \notag \\
& + 2C_3\big( 
  T^{\theta_3} \|u^N\|_{X^{-(1-a), a, T}}^3 
+ T^{\theta_4} \|\mathcal{N}^N\|_{X^{-a, a, T}}^3  \\
& +T^{\theta_5} \|\mathcal{N}^N\|_{X^{-a, a, T}}^2 
+T^{\theta_6} \|\mathcal{N}^N\|_{X^{-a, a, T}} 
+T^{\theta_7} \big). \notag
\end{align}

\noi
Moreover, for $N > M$, we have 
\begin{align} \label{GAMMA2}
\| u^N -  u^M &  \|_{X^{-(1-a), a, T}} 
   = \| \G^N u^N - \G^M u^M \|_{X^{-(1-a), a, T}} \notag \\
      \leq  & \  C_1 \|u_0^N - u_0^M\|_{H^{-(1-a)}}  \notag \\
 & + \tfrac{1}{2} C_2 T^{\theta_1} 
 \big(\|u^N\|_{X^{-(1-a), a, T}}+\|u^M\|_{X^{-(1-a), a, T}}\big)
 \|u^N-u^M\|_{X^{-(1-a), a, T}} \notag \\
& +   C_4 T^{\theta_3} \big(\|u^N\|_{X^{-(1-a), a, T}}^2
+\|u^M\|_{X^{-(1-a), a, T}}^2 \big) \|u^N-u^M\|_{X^{-(1-a), a, T}} \\
& + C_5 T^{\theta_4} \big(\|\mathcal{N}^N\|_{X^{-a, a, T}}^2  
+\|\mathcal{N}^M\|_{X^{-a, a, T}}^2\big)\|\mathcal{N}^N- \mathcal{N}^M\|_{X^{-a, a, T}}\notag\\
& +C_6 M^{-\eps_1}T^{\theta_5} \|\mathcal{N}^N - \mathcal{N}^M\|_{X^{-a, a, T}}^2 
 +C_6 M^{-\eps_1}T^{\theta_5} \|\mathcal{N}^N\|_{X^{-a, a, T}}^2 \notag \\
& +C_7 M^{-\eps_2}T^{\theta_6} \|\mathcal{N}^N- \mathcal{N}^M\|_{X^{-a, a, T}} 
+C_7 M^{-\eps_2}T^{\theta_6} \|\mathcal{N}^N\|_{X^{-a, a, T}} 
+C_8 M^{-\eps_3}T^{\theta_7}.  \notag
\end{align}

\noi
Also, we have
\begin{align} \label{GAMMA4}
\|\mathcal{N}^N -  & \mathcal{N}^M \|_{X^{-a, a, T}}\notag \\
 \leq & \  \tfrac{1}{2} C_2 T^{\theta_1} 
 \big(\|u^N\|_{X^{-(1-a), a, T}}+\|u^M\|_{X^{-(1-a), a, T}}\big)
 \|u^N-u^M\|_{X^{-(1-a), a, T}} \notag \\
 & + \ C_4 T^{\theta_3} \big(\|u^N\|_{X^{-(1-a), a, T}}^2
+\|u^M\|_{X^{-(1-a), a, T}}^2 \big) \|u^N-u^M\|_{X^{-(1-a), a, T}} \notag \\
& + C_5 T^{\theta_4} \big(\|\mathcal{N}^N\|_{X^{-a, a, T}}^2  
+\|\mathcal{N}^M\|_{X^{-a, a, T}}^2\big)\|\mathcal{N}^N- \mathcal{N}^M\|_{X^{-a, a, T}}\\
& +C_6 M^{-\eps_1}T^{\theta_5} \|\mathcal{N}^N - \mathcal{N}^M\|_{X^{-a, a, T}}^2 
 +C_6 M^{-\eps_1}T^{\theta_5} \|\mathcal{N}^N\|_{X^{-a, a, T}}^2 \notag \\
& +C_7 M^{-\eps_2}T^{\theta_6} \|\mathcal{N}^N- \mathcal{N}^M\|_{X^{-a, a, T}} 
+C_7 M^{-\eps_2}T^{\theta_6} \|\mathcal{N}^N\|_{X^{-a, a, T}} 
+C_8 M^{-\eps_3}T^{\theta_7}.  \notag
\end{align}

\noi
Note that in estimating the difference $\G^N u^N-\G^M u^M$ on $A_1$,  one needs to consider 
\begin{equation}\label{GAMMA3}
 \wt{I}_{-a, 1-a} := \| \mathcal{N}_1(u^N, u^N) - \mathcal{N}_1(u^M, u^M)\|_{-a, 1-a} 
\end{equation}

\noi
as in \cite{BO3}.
We can follow the argument on pp.135-136 in \cite{BO3}, yielding the third term in \eqref{GAMMA2}, 
except for $R_a$ defined in \eqref{262}. 
As for $R_a$, we can write 
\begin{align} \label{WTN1}
\mathcal{N}(\mathcal{N}(u,u), u)- \mathcal{N}(\mathcal{N}(v,v), v) 
= \mathcal{N}(\mathcal{N}(u+v,u-v),u)
+ \mathcal{N}(\mathcal{N}(v,v), u - v)
\end{align}

\noi as in (3.4) in \cite{BO3},
and then we can repeat the computation done for $R_a$, 
yielding the last six terms in \eqref{GAMMA2}.

Recall the large deviation estimate:
$ P(\|u_0(\omega)\|_{H^{s}} >K ) <e^{-c K^2}$ 
for $s < -\frac{1}{2}$ and sufficiently large $K$.
Given small $T>0$, let $K = ( 2C_1 C_3 T^{\theta_2})^{-1}$.
Then, 
defining $\Omega_T^{(2)}$ by \[\Omega_T^{(2)} := \{ \omega: \|u_0(\omega)\|_{H^{-(1-a)}} \leq K\} \]

\noi
we have $P\big((\Omega_T^{(2)})^c\big) < e^{-\frac{c}{T^\beta}}$
for some $c, \beta > 0$.
Moreover, by letting $R= 2C_1 \|u_0\|_{H^{-(1-a)}} $,
we have
$C_3 T^{\theta_2}R \leq \frac{1}{2}$ on $\Omega_T^{(2)}$.
Finally, let $\Omega_T = \Omega_T^{(1)}\cap \Omega_T^{(2)}$.
Then, by choosing $T$ sufficiently small, 
we see that for $\omega \in \Omega_T$, 
smooth global solutions $u^N(\omega) $ (and $\mathcal{N}^N$) with initial data $u^N_0(\omega)$ converge
in $X^{-(1-a), a, T}$
(in $X^{-a, a, T}$, respectively.)
For example, 
if we choose $T$ by 
\begin{align*} 
    T = \inf \big\{ t> 0: \tfrac{1}{2} C_2 t^{\theta_1} R 
  + 2C_3(  t^{\theta_3} R^3 
+ t^{\theta_4} R^3  & +t^{\theta_5} R^2 +t^{\theta_6} R 
+t^{\theta_7} )  \geq \tfrac{1}{2}R \big\},
 \end{align*}

\noi
then 
\eqref{GAMMA1} and \eqref{GAMMA5}
along with continuity argument
show that $\|u\|_{X^{-(1-a), a, T}}\leq R$
and 
$\|\mathcal{N}(u, u) \|_{X^{-a, a, T}} \leq R$.
From \eqref{GAMMA2} and \eqref{GAMMA4},
one obtains a different condition 
for $T$.
We point out that 
the nonlinear part $\mathcal{N}_j(u^N, u^N)$, $j = 1, 2$,  
converges in a stronger space $X^{-a, 1-a, T}$.
See \eqref{Isb} and \eqref{WT4}.

Let $u$ denote the limit. We still need to show, for $\omega \in \Omega_T$, 
\begin{itemize}
\item[(i)] $u$ is indeed a solution to \eqref{KDV} with $u_0 \in H^{-(1-a)}(\T)$
given by \eqref{Q IV}.
\item[(ii)] $u \in C([-T, T]; H^{-(1-a)})$.
\item[(iii)]  uniqueness of solutions.
\end{itemize}

\noi
The argument for (i) and (ii) 
exactly follows the corresponding argument in \cite{OH6},
and thus we omit details. 
It follows from \eqref{WT0} with $b = \frac{1}{2}+\dl$, 
\eqref{WT4}, and symmetry between $\s_1$ and $\s_2$
that 
\[
\begin{cases}
S(t) u_0 \in X^{-(1-a), \frac{1}{2}+\dl, T}
\subset C([-T, T]; H^{-(1-a)})\\
\mathcal{N}_1(u, u) +
\mathcal{N}_2(u, u) \in   X^{-a, \frac{1}{2}+\dl, T} 
\subset C([-T, T]; H^{-a}).
\end{cases}\]

\noi
As for $\mathcal{N}_0(u, u)$, i.e. $\s_0 = \MAX$,
we can repeat the argument in \cite{OH6}
by separately estimating
the contributions
on $A :=\{ \max(\s_1, \s_2) \gtrsim \jb{n n_1 n_2}^\frac{1}{100}\}$
and $A^c$.

Now, let us discuss uniqueness.
In the following, fix $\al \in ( \al_0, 0]$, where $\al_0$ is as in \eqref{al}.
Consider a ``nice'' initial condition $u_0^*: = u_0(\omega^*)$
for some $\omega^* \in \Omega_T$
so that 
a solution $u^*$ exists on $[-T, T]$ for the initial condition $u_0^*$
(along with the estimate in Subsection 4.2.)
Then, for some $\eps, \beta > 0$, we have
$\sup_{n\ne 0} \jb{n}^{-\eps} |g_n(\omega^*)| \leq C T^{-\frac{\beta}{2}}$. 
Now, let
\[A_{\g, T}= \{v_0: 
\sup_n \jb{n}^{-\eps}\big|\ft{v}_0(n) -  |n|^{-\al}g_n(\omega^*)\big| \leq \g C T^{-\frac{\beta}{2}}\}.\]

\noi
Then, we have 
$\sup_n \jb{n}^{-\eps} |\ft{v}_0(n)| \leq (1+\g) C T^{-\frac{\beta}{2}}$
on $A_{\g, T}$.
Moreover, for $v_0 \in A_{\g, T}$, we have
\begin{align*}
\| v_0 - u_0^*\|_{H^{-(1-a)}}
& = \bigg(\sum_{n\ne0} \jb{n}^{-1-2\dl}\big|\ft{v}_0(n) - |n|^{-\al} g_n (\omega^*)\big|^2\bigg)^\frac{1}{2}\\
& \lesssim \sup_n \jb{n}^{-\eps}\big|\ft{v}_0(n) -  |n|^{-\al} g_n(\omega^*)\big|
\leq \g C T^{-\frac{\beta}{2}}
\leq \tfrac{1}{10}\|u_0^*\|_{H^{-(1-a)}}
\end{align*}

\noi
by choosing $\eps < \dl$ and $\g$ sufficiently small.
Hence, proceeding as before, 
we can construct solutions $v$ with initial data $v_0\in A_{\g, T}$,
satisfying \eqref{GAMMA1} and \eqref{GAMMA5}
(after making a time interval slightly shorter.)
Consider the difference between $v$ and $u^* = u(\omega^*)$.
From a slight modification of the argument in Subsection 4.2,  
we have
\begin{align*} 
\| v -  u^* &  \|_{X^{-(1-a), a, T}} 
    \notag \\
 \leq  & \  C_1 \|v_0 - u_0^*\|_{H^{-(1-a)}}  \notag \\
 & + \tfrac{1}{2} C_2 T^{\theta_1} 
 \big(\|v\|_{X^{-(1-a), a, T}}+\|u^*\|_{X^{-(1-a), a, T}}\big)
 \|v-u^*\|_{X^{-(1-a), a, T}} \notag \\
& +   C_4 T^{\theta_3} \big(\|v\|_{X^{-(1-a), a, T}}^2
+\|u^*\|_{X^{-(1-a), a, T}}^2 \big) \|v-u^*\|_{X^{-(1-a), a, T}} \\
& + C_5 T^{\theta_4} \big(\|\mathcal{N}\|_{X^{-a, a, T}}^2  
+\|\mathcal{N}^*\|_{X^{-a, a, T}}^2\big)\|\mathcal{N}- \mathcal{N}^*\|_{X^{-a, a, T}}\notag\\
& +C_6 T^{\theta_5} \|\mathcal{N} - \mathcal{N}^*\|_{X^{-a, a, T}}^2 
 +C_6  T^{\theta_5} \|\mathcal{N}^*\|_{X^{-a, a, T}}^2 \sup_n \jb{n}^{-\eps}\big|\ft{v}_0(n) - |n|^{-\al} g_n(\omega^*)\big|\notag \\
& +C_7 T^{\theta_6} \|\mathcal{N}- \mathcal{N}^*\|_{X^{-a, a, T}} 
+C_7  T^{\theta_6} \|\mathcal{N}^*\|_{X^{-a, a, T}} \sup_n \jb{n}^{-\eps} \big|\ft{v}_0(n) -  |n|^{-\al}g_n(\omega^*)\big|\\
& +C_8  T^{\theta_7} \sup_n \jb{n}^{-\eps}\big|\ft{v}_0(n) - |n|^{-\al} g_n(\omega^*)\big|, \notag
\end{align*}

\noi
where
$\mathcal{N} = \mathcal{N}(v, v)$
and $\mathcal{N}^* = \mathcal{N}(u^*, u^*)$.
A similar estimate holds for
$\|\mathcal{N}^N -   \mathcal{N}^M \|_{X^{-a, a, T}}$.
As a consequence, we obtain
\[\| v -  u^*   \|_{X^{-(1-a), a, T}} \lesssim C(T, R) \sup_n \jb{n}^{-\eps}
\big|\ft{v}_0(n) -  |n|^{-\al}g_n(\omega^*)\big| = o(1)\] 

\noi
as $\g \to 0$.
Note that $u_0^* \in A_{\g, T}$ for any $\g > 0$.
Hence, $u^*$ is a unique solution.
One can similarly obtain 
\[\| v(t) -  u^* (t)  \|_{C([-T, T]; H^{-(1-a)})} \lesssim C(T, R) \sup_n \jb{n}^{-\eps}
\big|\ft{v}_0(n) -  |n|^{-\al} g_n(\omega^*)\big|.\] 

\noi
This provides a weak form of continuous dependence.
Lastly, when $\al = 0$,  
$u_0$ in \eqref{IVKDV} corresponds to the mean-zero Gaussian white noise.
In this case, one can extend local-in-time solutions
to global ones by invariance of the (finite dimensional) white noise.
See \cite{BO7,  OH7} for this part of discussion.

\section{On the Szeg\"o equation}

In this section, we consider the {\it dispersionless}
cubic Szeg\"o equation \eqref{Szego}
and present the proof of Proposition \ref{PROP:Szego}.
In particular, we show that unlike \cite{BO7} and \cite{CO1}, there is no gain 
of regularity
even if 
we take initial data to be random of the form \eqref{IVSZ}.
%

If we were to proceed as in \cite{BO7} and \cite{CO1},
we would need to estimate the ${C([-T, T]; H^s_+)}$-norm of $\mathcal{N}(u_1, u_2, u_3)$ 
defined in \eqref{SzDuhamel} for some $s\geq \frac{1}{2}$, 
assuming  $u_j$ is either of  type
\begin{itemize}
\item[(I)] random, less regular:
\[u_j (x, t) = \eta_{_T}(t)  u_0^\omega
= \eta_{_T}(t) \sum_{n \geq0 } \frac{ g_n(\omega)}{\sqrt{1+|n|^{2\al}}} e^{inx}, \ \ \text{or}
\]
\item[(II)] deterministic, smoother:\footnote{We could consider a different norm for type (II).
However, it is not relevant for the following discussion.}
\[u_j = v_j \text{ with } \|v_j \|_{C([-T, T]; H^s_+)} \leq 1,\]
\end{itemize}

In view of the well-posedness result in $H^{\frac{1}{2}}_+(\T)$,
we consider $\al \leq 1$ in the following.
(See \eqref{IVSZ}.)
We  show that the contribution from all type (I) is infinite a.s.
for $\al \leq 1$.
For notational simplicity, we use $\jb{n}^\al$ for $\sqrt{1+|n|^{2\al}}$.
Note that all the summations take place over non-negative indices,
i.e. $n\geq 0 $ and $ n_j \geq 0$, $j = 1, 2, 3$.

Now, assume $u_j$ is of type (I), $j = 1, 2, 3$.
Then, by separating the spatial and temporal components, we have
\begin{align}
\big\| \mathcal{N}(u_1,  u_2, u_3) & \big\|^2_{C([-T, T]; H^s_+)} 
 = C_T \bigg\| \jb{n}^s \Big|\sum_{n = n_1 - n_2 + n_3}
 \frac{g_{n_1}}{\jb{n_1}^\al}\frac{\cj{g_{n_2}}}{\jb{n_2}^\al}\frac{g_{n_3}}{\jb{n_3}^\al}
 \Big| \bigg\|_{l^2_n(\Z_{\geq 0})}^2\notag \\
 & = C_T \sum_n \jb{n}^{2s} 
 \sum_{\substack{n = n_1 - n_2 + n_3\\ \hphantom{nn} = m_1 - m_2 + m_3}}
\frac{g_{n_1}}{\jb{n_1}^\al}\frac{\cj{g_{n_2}}}{\jb{n_2}^\al}\frac{g_{n_3}}{\jb{n_3}^\al}
\frac{\cj{g_{m_1}}}{\jb{m_1}^\al}\frac{g_{m_2}}{\jb{m_2}^\al}\frac{\cj{g_{m_3}}}{\jb{m_3}^\al}.
\label{HHsum}
\end{align}

\noi
In the following, we first prove 
Proposition \ref{PROP:Szego} (a), 
showing that  \eqref{HHsum} is infinite a.s.~for $\al \in (\frac{1}{2}, 1]$
and $s \in [\al - \frac{1}{2}, \al - \frac{1}{4}]$.

We say that we have a {\it pair} 
if we have $n_j = m_j$ for some $j = 1, 2, 3$,
(or if we have $n_1 = m_3$ or $n_3 = m_1$.)
Then, we can separate the sum in \eqref{HHsum}
into three cases: (a) 3 pairs, (b) 1 pair, (c) no pair.
We estimate the contribution from each case in the following.

\medskip
\noi
$\bullet$ {\bf Case (a):} 3 pairs.
\quad
In this case, the contribution to \eqref{HHsum} is bounded from below by
\begin{align} \label{HHsuma}
\sum_n \jb{n}^{2s} 
 \sum_{n = n_1 - n_2 + n_3}
\frac{|g_{n_1}|^2}{\jb{n_1}^{2\al}}\frac{|g_{n_2}|^2}{\jb{n_2}^{2\al}}\frac{|g_{n_3}|^2}{\jb{n_3}^{2\al}}.
\end{align}

\noi
Now, consider the contribution from $n = n_1$ and  $n_2 = n_3$.
For $\al > \frac{1}{4}$, we have
\[c_\omega := \sum_{n_2} \jb{n_2}^{-4\al} |g_{n_2}(\omega)|^4 < \infty, \  \text{ a.s.}\]

\noi
since 
$\mathbb{E}\big[ \sum_{n_2} \jb{n_2}^{-4\al} |g_{n_2}|^4\big]
\sim \sum_{n_2} \jb{n_2}^{-4\al} <\infty$.
Note that $c_\omega > 0$ a.s.
Let $F_j(\omega):= 2^{-j} \sum_{|n|\sim 2^j}  |g_n(\omega)|^2$.
Then, $F_j(\omega)$ 
converges to a positive constant a.s. by strong law of large numbers.
Hence, 
for $\al \leq s +\frac{1}{2}$, we have
\[
\sum_{|n|\sim 2^j} \jb{n}^{2s-2\al} |g_n(\omega)|^2
\sim 2^{j(2s-2\al+1)} F_j(\omega) \not\to 0, \ \text{ a.s.}\]

\noi
Therefore, we have
\begin{align*}
\eqref{HHsuma} \geq c_\omega \sum_{n} \jb{n}^{2s-2\al} |g_n(\omega)|^2
= \infty, \ \text{ a.s.}
\end{align*}

\noi
for $\al \leq s +\frac{1}{2}$.
In particular, when $s =\frac{1}{2}$, the contribution from this case is divergent for $\al \leq 1$.
This already shows that there is no nonlinear smoothing
even if we consider random initial data of the form \eqref{IVSZ}.

Suppose that we have one pair $n_1 = m_1$.
Moreover, assume $n_2 = n_3$.
Then, we have $n = n_1 = m_1$
and thus $m_2 = m_3$. 
Proceeding in a similar manner as above, 
we see that the contribution to \eqref{HHsum} is given by
\begin{align} 
\sum_n \jb{n}^{2s} \frac{|g_{n}|^2}{\jb{n}^{2\al}}
 \sum_{n_2, m_2}
\frac{|g_{n_2}|^2}{\jb{n_2}^{2\al}}
\frac{|g_{m_2}|^2}{\jb{m_2}^{2\al}} = \infty, \ \text{ a.s.}
\end{align}

\noi
for $\al \leq s+\frac{1}{2}$.

For completeness of the argument, 
we give a brief discussion to show that 
the contributions from other cases are finite (at least for $\al >\frac{1}{2}$
and $s \leq \al - \frac{1}{4}$.)

\medskip
\noi
$\bullet$ {\bf Case (b):} 1 pair.
\quad
We only discuss the case when $n_2 = m_2$ and $\{n_1, n_3\}\ne \{m_1, m_3\}$.
Other cases follow in a similar manner
(except for the case discussed above.)
In this case, the contribution to \eqref{HHsum} is given by
\begin{align*}
R_1 (\omega): =  \sum_n \jb{n}^{2s} 
 \sum_{n_2} \frac{|g_{n_2}|^2}{\jb{n_2}^{2\al}}
 \sum_{\substack{n + n_2 =n_1+ n_3\\ \hphantom{nnnnll} = m_1  + m_3}}
\frac{g_{n_1}}{\jb{n_1}^\al}\frac{g_{n_3}}{\jb{n_3}^\al}
\frac{\cj{g_{m_1}}}{\jb{m_1}^\al}\frac{\cj{g_{m_3}}}{\jb{m_3}^\al}.
\end{align*}

\noi
By computing the second moment, we have
\begin{align}
\label{HHsumb2}
\mathbb{E}\big[|R_1|^2 \big] 
 =  \mathbb{E}\bigg[ & \sum_n \jb{n}^{2s} 
 \sum_{n_2} \frac{|g_{n_2}|^2}{\jb{n_2}^{2\al}}
  \sum_{\substack{n + n_2 =n_1+ n_3\\ \hphantom{nnnnll}= m_1  + m_3}}
\frac{g_{n_1}}{\jb{n_1}^\al}\frac{g_{n_3}}{\jb{n_3}^\al}
\frac{\cj{g_{m_1}}}{\jb{m_1}^\al}\frac{\cj{g_{m_3}}}{\jb{m_3}^\al} \notag \\
\times
& \sum_{\wt{n}} \jb{\wt{n}}^{2s} 
 \sum_{\wt{n}_2} \frac{|g_{\wt{n}_2}|^2}{\jb{\wt{n}_2}^{2\al}}
  \sum_{\substack{\wt{n} + \wt{n}_2 =\wt{n}_1+ \wt{n}_3\\ 
  \hphantom{nnnnll}= \wt{m}_1  + \wt{m}_3}}
\frac{\cj{g_{\wt{n}_1}}}{\jb{\wt{n}_1}^\al}\frac{\cj{g_{\wt{n}_3}}}{\jb{\wt{n}_3}^\al}
\frac{g_{\wt{m}_1}}{\jb{\wt{m}_1}^\al}\frac{g_{\wt{m}_3}}{\jb{\wt{m}_3}^\al}\bigg].
  \end{align}
  
\noi
We have nontrivial contribution in \eqref{HHsumb2} 
when $(n_1, n_3, \wt{m}_1, \wt{m}_3)
= (\wt{n}_1, \wt{n}_3, m_1, m_3)$
(up to permutation.)
By assumption, we have $\{n_1, n_3\}\ne \{m_1, m_3\}$
and $\{\wt{n}_1, \wt{n}_3\}\ne \{\wt{m}_1, \wt{m}_3\}$.
Then, it follows that $(n_1, n_3) = (\wt{n}_1, \wt{n}_3)$
and $(m_1, m_3) = ( \wt{m}_1, \wt{m}_3)$
(up to permutation.)
Hence, we have
\begin{align} \label{HHsumb3}
\eqref{HHsumb2}
\sim 
 \sum_n \jb{n}^{2s} 
 & \sum_{n_2} \frac{1}{\jb{n_2}^{2\al}}
  \sum_{ n+ n_2 = \wt{n}+\wt{n}_2  } \jb{\wt{n}}^{2s} 
 \frac{1}{\jb{\wt{n}_2}^{2\al}} \notag \\
&  \times
 \sum_{\substack{n + n_2 =n_1+ n_3\\ \hphantom{nnnnll}= m_1  + m_3}}
\frac{1}{\jb{n_1}^{2\al}}\frac{1}{\jb{n_3}^{2\al}}
\frac{1}{\jb{m_1}^{2\al}}\frac{1}{\jb{m_3}^{2\al}}.
  \end{align}

\noi
Without loss of generality, assume $n_1 \gtrsim n$, since $n_1+ n_3 = n+n_2 \geq n$.
Then, for fixed $n$ and $n_2$, we have
\begin{align*}
\jb{n}^{2s}
\sum_{n + n_2 =n_1+ n_3}
\frac{1}{\jb{n_1}^{2\al}}\frac{1}{\jb{n_3}^{2\al}}
\lesssim \sum_{n + n_2 =n_1+ n_3}
\frac{1}{\jb{n_1}^{2(\al-s)}}\frac{1}{\jb{n+n_2 -n_1}^{2\al}}
\lesssim \frac{1}{\jb{n+n_2}^{\frac{1}{2}+}}
\end{align*}

\noi
for $2(\al - s) \geq \frac{1}{2}$ and $\al > \frac{1}{2}$, 
where the last inequality follows
from a slight modification of \cite[Lemma 2.2]{TZ3}.
The same inequality holds when $n$ and $n_j$ are replaced by $\wt{n}$ and $\wt{n}_j$.
Hence, we have
\begin{align}
\eqref{HHsumb3}
\lesssim 
\sum_n \jb{n}^{-1-} 
 \sum_{n_2} \jb{n_2}^{-2\al} \sum_{\wt{n}_2} \jb{\wt{n}_2}^{-2\al}<\infty.
 \end{align}

\noi
Therefore, we have $R_1(\omega) < \infty$ a.s.
for $\al >\frac{1}{2}$
and $s \leq \al - \frac{1}{4}$.

\medskip
\noi
$\bullet$ {\bf Case (c):} no pair.
\quad
In this case, the contribution to \eqref{HHsum} is given by
\begin{align}\label{HHsumc0}
R_2(\omega) := \sum_n \jb{n}^{2s} 
 \sum_{*}
\frac{g_{n_1}}{\jb{n_1}^\al}\frac{\cj{g_{n_2}}}{\jb{n_2}^\al}\frac{g_{n_3}}{\jb{n_3}^\al}
\frac{\cj{g_{m_1}}}{\jb{m_1}^\al}\frac{g_{m_2}}{\jb{m_2}^\al}\frac{\cj{g_{m_3}}}{\jb{m_3}^\al},
\end{align}

\noi
where
$* = \{ n = n_1 - n_2 + n_3 = m_1 - m_2 + m_3, \text{ no pair}\}$.
First, suppose $n = n_1$,
and thus $n_2 = n_3$.
Then, the contribution to \eqref{HHsumc0} is given by
\begin{align*}
 \sum_n \jb{n}^{2s} \sum_{n_2} \frac{|g_{n_2}|^2}{\jb{n_2}^{2\al}}
 \sum_{n = m_1 - m_2 + m_3}
\frac{g_{n}}{\jb{n}^\al}
\frac{\cj{g_{m_1}}}{\jb{m_1}^\al}
\frac{g_{m_2}}{\jb{m_2}^\al}
\frac{\cj{g_{m_3}}}{\jb{m_3}^\al}.
\end{align*}

\noi
By computing the second moment as before, we have
\begin{align*}
\mathbb{E}\big[|R_2|^2 \big] 
\sim \sum_n \jb{n}^{2s}  \max(\jb{n}^{2s},\jb{m_2}^{2s})
& \sum_{n_2}\frac{1}{\jb{n_2}^{2\al}}
\sum_{\wt{n}_2}\frac{1}{\jb{\wt{n}_2}^{2\al}}\\
& \times \sum_{n = m_1 - m_2 + m_3}
 \frac{1}{\jb{n}^{2\al}}
\frac{1}{\jb{m_1}^{2\al}}
\frac{1}{\jb{m_2}^{2\al}}
\frac{1}{\jb{m_3}^{2\al}}.
\end{align*}

\noi
Now, we can follow the argument in Case (b) to show that $\mathbb{E}\big[|R_2|^2 \big] <\infty$.

In the following, assume that $n_1, n_3, m_1, m_3 \ne n$.
In this case, we have
\begin{align}
\label{HHsumc}
\mathbb{E}\big[|R_2|^2 \big] 
&  =  \mathbb{E}\bigg[  \sum_n \jb{n}^{2s} 
   \sum_{*}
\frac{g_{n_1}}{\jb{n_1}^\al}\frac{\cj{g_{n_2}}}{\jb{n_2}^\al}\frac{g_{n_3}}{\jb{n_3}^\al}
\frac{\cj{g_{m_1}}}{\jb{m_1}^\al}\frac{g_{m_2}}{\jb{m_2}^\al}\frac{\cj{g_{m_3}}}{\jb{m_3}^\al} \notag \\
& \hphantom{XXX}
\times \sum_{\wt{n}} \jb{\wt{n}}^{2s} 
   \sum_{\wt{*}}
\frac{\cj{g_{\wt{n}_1}}}{\jb{\wt{n}_1}^\al}
\frac{g_{\wt{n}_2}}{\jb{\wt{n}_2}^\al}\frac{\cj{g_{\wt{n}_3}}}{\jb{\wt{n}_3}^\al}
\frac{g_{\wt{m}_1}}{\jb{\wt{m}_1}^\al}\frac{\cj{g_{\wt{m}_2}}}{\jb{\wt{m}_2}^\al}
\frac{g_{\wt{m}_3}}{\jb{\wt{m}_3}^\al}\bigg]\\
& \sim \sum_n \jb{n}^{2s} 
   \sum_{*}
\jb{\wt{n}}^{2s}
\prod_{j = 1}^3
\frac{1}{\jb{n_j}^{2\al}}
\frac{1}{\jb{m_j}^{2\al}}\notag 
  \end{align}

\noi
where
$\wt{*} = \{ \wt{n} = \wt{n}_1 - \wt{n}_2 + \wt{n}_3 = \wt{m}_1 - \wt{m}_2 + \wt{m}_3, \text{ no pair },
 \wt{n}_1, \wt{n}_3, \wt{m}_1, \wt{m}_3\ne \wt{n} \}$.
Note that there is no summation for $\wt{n}$ 
since it is determined by the values $n_j$ and $m_j$.
Without loss of generality, assume $n_1 \gtrsim n$.
Then, by (a slight modification of) \cite[Lemma 2.2]{TZ3},
we have
\begin{align*}
 \jb{n}^{2s} \sum_* \prod_{j = 1}^3
\frac{1}{\jb{n_j}^{2\al}}
& \lesssim \sum_{n_2} \frac{1}{\jb{n_2}^{2\al}}
\sum_{n_1} \frac{1}{\jb{n_1}^{2(\al-s)}}
\frac{1}{\jb{n+n_2-n_1}^{2\al}}\\
& \lesssim \sum_{n_2} \frac{1}{\jb{n_2}^{2\al}}
 \frac{1}{\jb{n+n_2}^{\frac{1}{2}+}}
\lesssim 
\frac{1}{\jb{n}^{\frac{1}{2}+}}
\end{align*}

\noi
for $2(\al - s) \geq \frac{1}{2}$ and $\al >\frac{1}{2}$.
Hence, we have $\mathbb{E}\big[|R_2|^2 \big] <\infty$
and thus $R_2(\omega) < \infty$ a.s.
for $\al >\frac{1}{2}$
and $s \leq \al - \frac{1}{4}$.

In general, i.e. if $s > \al - \frac{1}{4}$,
then we can choose $\wt{s} < s $ such that 
$\wt{s}  \in [\al - \frac{1}{2}, \al - \frac{1}{4}]$.
Then, from the previous computation with $\wt{s}$ instead of $s$, 
we have
\begin{align*}
\big\| \mathcal{N}(u_1,  u_2, u_3)  \big\|_{C([-T, T]; H^s_+)}
\geq \big\| \mathcal{N}(u_1,  u_2, u_3)  \big\|_{C([-T, T]; H^{\wt{s}}_+)}
= \infty, \quad \text{a.s.}
\end{align*}

\noi
This proves Part (a) of Proposition \ref{PROP:Szego}.

\medskip
Part (b) of of Proposition \ref{PROP:Szego} follows easily 
by taking an expectation of \eqref{HHsum}.
After taking an expectation, only  the case with 3 pairs remains.
Assume $n_1 = \max(n_1, n_2, n_3)$ in the following. 
Then, we have $n_1 \gtrsim n$.
With $s = \al - \frac{1}{2}-\eps$, we have
\begin{align*}
\mathbb{E} \Big[\big\| \mathcal{N}(u_1,  u_2, u_3)  \big\|^2_{C([-T, T]; H^s_+)} \Big]
& \lesssim 
\sum_n \jb{n}^{2s} 
 \sum_{n = n_1 - n_2 + n_3}
\frac{1}{\jb{n_1}^{2\al}}\frac{1}{\jb{n_2}^{2\al}}\frac{1}{\jb{n_3}^{2\al}}\\
& \lesssim 
\sum_n 
 \sum_{n = n_1 - n_2 + n_3}
\frac{1}{\jb{n_1}^{1+\eps}}\frac{1}{\jb{n_2}^{2\al}}\frac{1}{\jb{n_3}^{2\al}}
<\infty,
\end{align*}

\noi
as long as $\al > \frac{1}{2}$.

\begin{remark}\label{REM:last}\rm
In studying nonlinear smoothing under randomization for
the Wick ordered cubic NLS in \cite{CO1}, 
we needed to control a term similar to \eqref{HHsum}
when we estimate the contribution from all type (I) to the nonlinearity in the $X^{s, -\frac{1}{2}+}$ norm.
On the one hand, if $\s_0 := \jb{ \tau - n^2}$ is large, 
the estimate was trivial.
On the other hand, if $\s_0 := \jb{ \tau - n^2}$ is small, then
the relation: 
\[ \s_0 - \s_1 + \s_2 - \s_3
= -n^2 + n_1^2 - n_2^2 + n_3^2= -2(n_2 - n_1)(n_2-n_3),\]

\noi
where $\s_j := \jb{\tau_j - n_j^2}\lesssim 1$
with $n = n_1 - n_2 + n_3$ and $\tau = \tau = \tau_1 - \tau_2 + \tau_3$,
imposed a restriction on the summation.
See \cite{CO1} for details.
However, due to the lack of dispersion for \eqref{Szego},
there is no restriction on the summation \eqref{HHsum}.

In the case of KdV, 
although a direct estimate on the integral formulation failed
to show any nonlinear smoothing (see Subsection 4.1),
we could show that there is a gain of regularity
by considering the second iteration.
This is due to the fact that
$\s_1 = \jb{\tau_1 -n_1^3}\gtrsim \jb{n n_1 n_2}$ appears in the denominator
in the second iteration.
See \eqref{eq:I2}.
However, even if we consider the second iteration
for the cubic Szeg\"o equation \eqref{Szego}, 
it seems that we do not have any gain 
due to the lack of dispersion.
\end{remark}

\end{document}